\documentclass[12pt,fleqn,leqno]{article}
\usepackage{graphicx,color,enumerate,amssymb}
\usepackage{amsmath}
\usepackage{amssymb}
\usepackage{epsfig}
\usepackage{array}
\usepackage{float}

\makeatletter
\@addtoreset{equation}{section}
\renewcommand\thetable{\@arabic\c@table}
\renewcommand\thefigure{\@arabic\c@figure}
\long\def\@makecaption#1#2{%
  \vskip\abovecaptionskip
  \begin{center}%
  \sbox\@tempboxa{#1: #2}%
  \ifdim \wd\@tempboxa >\hsize
    #1: #2\par
  \else
    \global \@minipagefalse
    \hb@xt@\hsize{\hfil\box\@tempboxa\hfil}%
  \fi
  \end{center}%
  \vskip\belowcaptionskip}
\def\N{{\rm I\kern-.15em N}}
\def\R{{\rm I\kern-.2em R}}
\def\Z{{\rm Z\kern-.26em Z}}
\makeatother
\newtheorem{thm}{Theorem}[section]

\newtheorem{rem}[thm]{Remark}
\newtheorem{cor}[thm]{Corollary}
\newtheorem{prop}[thm]{Proposition}

\newcommand{\be}{\begin{eqnarray}}
\newcommand{\ee}{\end{eqnarray}}
\newcommand{\bq}{\begin{eqnarray*}}
\newcommand{\eq}{\end{eqnarray*}}

\newcommand{\bewend}{\hspace*{2mm}\rule{3mm}{3mm}}

\newcommand{\RR}{\mathbb{R}}
\newcommand{\PP}{\mathbb{P}}

\newcommand{\BE}{\mathbb{E}}
\newcommand{\e}{\textrm{e}}
\newcommand{\LL}{{\textrm L}^2}
\newcommand{\vertk}{\stackrel{{\cal D}}{\longrightarrow}}

\begin{document}
%%%%%%%%%%%%%%%%%%%%%%%%%%%%%%%%%%%%%%%%%%%%%%%%%%%%%%%%%%%%%%%%%%%%%%%%%%%
\begin{center}
{ \LARGE\sc Testing for normality in any dimension based on a partial differential equation involving the moment generating function\\
} \vspace*{0.5cm} {\large\sc Norbert Henze}$^{1}$, {\large\sc Jaco Visagie}$^{2}$\\ \vspace*{0.5cm}
{\it $^{1}$Institute of Stochastics, Karlsruhe Institute of Technology (KIT), Karlsruhe, Germany \\
henze@kit.edu} \\
{\it $^{2}$Department of Statistics, University of Pretoria, South Africa \\
jaco.visagie@up.ac.za}\\
%\today
\end{center}

{\small {\bf Abstract.}{ We use a system of first-order partial differential equations that characterize the moment generating function of the
$d$-variate standard normal distribution to construct a class of affine invariant tests for normality in any dimension.
We derive the limit null distribution of the resulting test statistics, and we prove consistency of the tests against general alternatives.
In the case $d>1$, a certain limit of these tests is connected with two measures of multivariate skewness.
The new tests show strong power performance when compared to
well-known competitors, especially against heavy-tailed distributions, and they are illustrated by means of a real data set.}\\
{\small {\it Keywords.}
Moment generating function; test for multivariate normality; direct sum of Hilbert spaces; multivariate skewness; weighted $L^2$-statistic \\
\vspace*{0.3cm} {2010 AMS classification:} 62H15, 62G20.
%\\ \\
%\newpage

%
%
%
%

%%%%%%%%%%%%%%%%%%%%%%%%%%%%%%%%%%%%%%%%%%%%%%%%%%%%%%%%%%%%%%%%%%%%%%%%%%%
\section{Introduction}\label{sec_1}

It is often of interest to check whether an observed $d$-dimensional dataset is compatible with the assumption of coming from a multivariate normal distribution. 
Such a model check is of practical interest to researchers, as many multivariate techniques rely on the assumption of multivariate normality (for short: multinormality).
As a consequence, it is not surprising that there is ongoing interest in testing for multinormality,
as can be witnessed by various recent papers on the subject, see for example \cite{MM2005}, \cite{SR2005}, \cite{FSN2007} and \cite{JV2014} as well as the references therein.
Research into the practical implementation of these tests is also undertaken regularly, see, e.g.,  \cite{KGZ2014} and \cite{JV2014} regarding
the implementation in the statistical software package R.

The purpose of this paper is not to review the multitude of tests that hitherto has been proposed (for an account of the state of the art regarding affine invariant procedures before 2002, see \cite{Hen2002}),
but to introduce and study a new class of tests that is based on a certain partial differential equation. To be specific, let
$X_1,X_2, \ldots, X_n, \ldots $ be independent and identically distributed (i.i.d.) copies of a $d$-dimensional random (column) vector $X$,  the distribution of which
is assumed to be absolutely continuous with respect to $d$-dimensional Lebesgue measure. All random vectors are defined on a common probability space $(\Omega,{\cal A},\mathbb{P})$.

Writing N$_d(\mu,\Sigma)$ for the $d$-dimensional normal distribution with mean vector $\mu$ and non-degenerate covariance matrix $\Sigma$
and ${\cal N}_d$ for the class of all non-degenerate $d$-variate normal distributions, a classical problem is to test the null hypothesis
\begin{equation}\label{NULL} H_0: \ {\rm{The \ law \ of}} \  X  \ {\rm{belongs\  to \ }} \ {\cal{N}}_d,
\end{equation}
against general alternatives. Since the class ${\cal N}_d$ is closed with respect to full rank affine transformations, any genuine test statistic $T_n = T_n(X_1,\ldots,X_n)$ based on $X_1,\ldots,X_n$ should also be invariant with respect to such transformations, i.e., we should have
$
T_n(AX_1+b, \ldots, AX_n+b)  =  T_n(X_1,\ldots,X_n)
$
for each nonsingular ($d \times d$)-matrix $A$ and each $b \in \mathbb{R}^d$, see \cite{Hen2002} for an account of the importance of affine invariance in connection
with testing for multivariate normality. Writing $\overline{X}_n = n^{-1}\sum_{j=1}^n X_j$ for the sample mean and
 $S_n = n^{-1}\sum_{j=1}^n (X_j - \overline{X}_n)(X_j - \overline{X}_n)^\top$ for the sample covariance matrix of $X_1,\ldots,X_n$,
 where the superscript $\top$ denotes the transpose of vectors and matrices, a necessary and sufficient condition for a test statistic $T_n$ to be affine invariant is that it is based on the scalar products
 \[
 Y_{n,i}^\top Y_{n,j} \ = \ (X_i - \overline{X}_n)^\top S_n^{-1} (X_j - \overline{X}_n), \qquad i,j \in \{1,\ldots, n\},
 \]
 where $Y_{n,j} = S_n^{-1/2}(X_j - \overline{X}_n)$, $j=1,\ldots, n$,
 are the so-called {\em scaled residuals} of $X_1,\ldots,X_n$, see \cite{Hen2002}.
Here, $S_n^{-1/2}$ denotes the unique symmetric square root of $S_n^{-1}$ which, due to the absolute continuity of the distribution of $X$, exists with probability one if $n \ge d+1$, see \cite{EP1973}. The latter condition is tacitly assumed to hold in what follows.

The novel idea for constructing a test of $H_0$ is the following: Suppose $X$ is a $d$-dimensional random vector, and assume that the moment generating function (MGF) $m(t) = \mathbb{E} [\exp(t^\top X)]$ exists for each
$t \in \mathbb{R}^d$ and satisfies the system of partial differential  equations
\begin{equation}\label{partial}
\frac{\partial m(t)}{\partial t_j}  = t_j m(t),  \quad t=(t_1,\ldots,t_d)^\top \in \mathbb{R}^d, \quad j=1,\ldots,d.
\end{equation}
Writing $y'(t)$ for the gradient of a function $y:\mathbb{R}^d \to \mathbb{R}$ at $t$, (\ref{partial}) may be succinctly written as
$
m'(t) = t m(t)$, $t \in \mathbb{R}^d$.
Taking into account that $m(0)=1$, the unique solution of this equation in the case $d=1$ is $m(t) = \exp(t^2/2)$, which is the MGF of the standard normal distribution.
If $d>1$, and we fix $t_2,\ldots,t_d$, the solution of (\ref{partial}) for $j=1$ is
\begin{equation}\label{partial2}
m(t) = c_2(t_2,\ldots,t_d) \cdot \textrm{e}^{t_1^2/2}
\end{equation}
for some function $c_2:\mathbb{R}^{d-1} \to \mathbb{R}$. We thus have $c_2(t_2,\ldots,t_d) = \textrm{e}^{- t_1^2/2} \cdot m(t)$, $t \in \mathbb{R}^d$,
which shows that $c_2$ is differentiable. Moreover,
\[
\frac{\partial}{\partial t_2} c_2(t_2,\ldots,t_d) = \textrm{e}^{- t_1^2/2} \cdot \frac{\partial}{\partial t_2} m(t).
\]
Using (\ref{partial}) with $j=2$ then gives
\[
t_2 = \frac{\frac{\partial}{\partial t_2} c_2(t_2,\ldots,t_d)}{c_2(t_2,\ldots,t_d)}, \quad t \in \mathbb{R}^d,
\]
the solution of which is $c_2(t_2,\ldots,t_d) = c_3(t_3,\ldots,t_d) \cdot \exp(t_2^2/2)$ for some function $c_3:\mathbb{R}^{d-2} \to \mathbb{R}$.
Inserting this expression into (\ref{partial2}) and continuing this way, we finally obtain
\[
m(t) = \prod_{j=1}^d \textrm{e}^{t_j^2/2} = \textrm{e}^{\|t\|^2/2}, \quad t \in \mathbb{R}^d,
\]
where $\| \cdot \|$ denotes   the Euclidean norm in $\mathbb{R}^d$. Notice that this unique solution of (\ref{partial}) is the MGF of the standard normal distribution N$_d(0,\textrm{I}_d)$, where $\textrm{I}_d$ is the unit matrix of order $d$.

If $X$ has some non-degenerate normal distribution, the scaled residuals $Y_{n,1},\ldots,Y_{n,n}$ should be approximately independent, with a distribution close to N$_d(0,\textrm{I}_d)$, at least for large $n$.
Hence, a natural approach for testing $H_0$ is to consider the empirical MGF
\begin{equation}\label{defmgf}
M_n(t) = \frac{1}{n} \sum_{j=1}^n \textrm{e}^{t^\top Y_{n,j}}
\end{equation}
of $Y_{n,1},\ldots,Y_{n,n}$, and to employ the
weighted $L^2$-statistic
\begin{equation}\label{teststat}
T_{n,\gamma} := n \int_{\mathbb{R}^d} \|M_n'(t) - t M_n(t)\|^2 \, w_\gamma(t)  \, \textrm{d}t,
\end{equation}
where
\begin{equation}\label{weightf}
w_\gamma(t) = \exp\left(- \gamma \|t\|^2\right), \qquad t \in \mathbb{R}^d.
\end{equation}
Rejection of $H_0$ is for large values of $T_{n,\gamma}$.
The role of $\gamma >0$ will be discussed later.

Using the relations
\begin{eqnarray*}
\int_{\mathbb{R}^d} \textrm{e}^{t^\top \alpha} \, w_\gamma(t) \, \textrm{d}t & = & \left(\frac{\pi}{\gamma}\right)^{d/2} \exp\left(\frac{\|\alpha\|^2}{4\gamma}\right),\\
\int_{\mathbb{R}^d} \textrm{e}^{t^\top \alpha} \, t^\top \alpha \, w_\gamma(t) \, \textrm{d}t & = & \left(\frac{\pi}{\gamma}\right)^{d/2} \, \frac{\|\alpha\|^2}{2\gamma} \, \exp\left(\frac{\|\alpha\|^2}{4\gamma}\right),\\
\int_{\mathbb{R}^d} \textrm{e}^{t^\top \alpha} \, \|t\|^2 \, w_\gamma(t) \, \textrm{d}t & = & \left(\frac{\pi}{\gamma}\right)^{d/2} \exp\left(\frac{\|\alpha\|^2}{4\gamma}\right) \left(\frac{d}{2 \gamma} + \frac{\|\alpha\|^2}{4\gamma^2}\right),
\end{eqnarray*}
and putting $Y_{n,j,k}^+ = Y_{n,j} + Y_{n,k}$,
the test statistic defined in (\ref{teststat}) takes the form
\begin{equation}\label{rep1}
T_{n,\gamma} \! = \! \frac{1}{n} \! \! \left(\frac{\pi}{\gamma}\right)^{d/2} \! \! \sum_{j,k=1}^n \! \! \! \exp \! \left(\! \frac{\|Y_{n,j,k}^+\|^2}{4 \gamma} \! \right) \! \! \! \left(\! Y_{n,j}^\top Y_{n,k} \! - \!  \frac{\|Y_{n,j,k}^+\|^2}{2\gamma}
\! + \!  \frac{d}{2\gamma} \! + \!  \frac{\|Y_{n,j,k}^+\|^2}{4\gamma^2} \! \right) \! \! ,
\end{equation}
which is amenable to computational purposes. Notice that $T_{n,\gamma}$ is affine invariant.

The remainder of this paper is structured as follows: Section \ref{seclimnull} deals with the convergence in distribution of $T_{n,\gamma}$ under $H_0$, and
Section \ref{secconsist} is devoted to the problem of consistency of the new test (which seems to be new even in the univariate case). In Section \ref{seclimgam}
we let $\gamma$ tend to infinity while keeping $n$ fixed.
Under this setting, $T_{n,\gamma}$ converges to a certain linear combination of two well-known measures of multivariate skewness.
In Section \ref{secsimul} we compare the finite-sample power behavior of the new test to that of several classical and recent tests for both univariate and multivariate normality.
Section \ref{secrealdata} illustrates the procedures by means of a real data set.
Section \ref{secconclus} presents some conclusions, while Section \ref{sectechnical} contains several technical proofs which do not form part of the main text.

%%%%%%%%%%%%%%%%%%%%%%%%%%%%%%%%%%%%%%%%%%%%%%%%%%%%%%%%%%%%%%%%%%%%%%%%%%%%%%%%%%%%%%%%%%%%%%%%%%%%%%%%%%%%%%%%%%%%%%%%%%%%%%%%%%%%%%%%%%%%%%%%%%%%%%%%%%%%%%%%%%%%%%%%%%%%%%%%%%%%%%%%%%%%%%%%%%%%%%%%%%%
%
%
%       The limit null distribution
%
%
%%%%%%%%%%%%%%%%%%%%%%%%%%%%%%%%%%%%%%%%%%%%%%%%%%%%%%%%%%%%%%%%%%%%%%%%%%%%%%%%%%%%%%%%%%%%%%%%%%%%%%%%%%%%%%%%%%%%%%%%%%%%%%%%%%%%%%%%%%%%%%%%%%%%%%%%%%%%%%%%%%%%%%%%%%%%%%%%%%%%%%%%%%%%%%%%%%%%%%%%%%%%

\section{The limit null distribution of $T_{n,\gamma}$}\label{seclimnull}
In this section, we derive the limit distribution of the test statistic $T_{n,\gamma}$ defined in (\ref{teststat}) under the null hypothesis (\ref{NULL}).
In view of affine invariance,
we will assume $\BE \left( X \right) = 0$ and $\BE \left( XX^\top \right) = \textrm{I}_d$. Put
\begin{equation}\label{defwn}
W_n(t) := \sqrt{n} \left(M_n'(t) - t M_n(t)\right) \ = \ \frac{1}{\sqrt{n}} \sum_{j=1}^n \e^{t^\top Y_{n,j}}\left(Y_{n,j} - t\right), \qquad t \in \RR^d,
\end{equation}
and let
$\textrm{L}^2 := \textrm{L}^2(\RR^d,{\cal B}^d,w_\gamma(t)\textrm{d}t)$ denote the separable Hilbert space of (equivalence classes of) measurable functions $f:\RR^d \rightarrow \RR$
that are square integrable with respect to the finite measure on the $\sigma$-field ${\cal B}^d$ of Borel sets of $\RR^d$, given by the weight function
$w_\gamma$. The inner product and the resulting norm on $\textrm{L}^2$ will be denoted by
\[
\langle f,g \rangle = \int_{\RR^d} f(t) \, g(t) \, w_\gamma(t) \, \textrm{d}t, \qquad \|f\|_{\textrm{L}^2} = \sqrt{\langle f,f \rangle},
\]
respectively. Since $W_n(t)$ in (\ref{defwn}) is $\RR^d$-valued, we consider the Hilbert space ${\cal H}$, which is the $d$-fold (orthogonal) direct sum ${\cal H} := \LL \oplus \cdots \oplus \LL$, consisting of
all ordered $d$-tuples $f= (f_1,\ldots,f_d) \in \LL \times \cdots \times \LL$, equipped with the  inner product
\[
\langle f, g \rangle_\oplus := \langle f_1,g_1\rangle_{\LL} + \ldots + \langle f_d,g_d\rangle_{\LL},
\]
where $f=(f_1,\ldots,f_d)$, $g = (g_1,\ldots, g_d) \in {\cal H}$. Notice that the norm $\|\cdot \|_\oplus$ on ${\cal H}$ satisfies
\[
\|f\|^2_\oplus = \sum_{i=1}^d \|f_i\|^2_{\LL} = \int_{\RR^d} \|f(t)\|^2\, w_\gamma(t) \, \textrm{d}t, \qquad f=(f_1,\ldots,f_d) \in {\cal H}.
\]
Since $W_n$ is a random element of ${\cal H}$ and $T_n = \|W_n\|^2_\oplus$, the aim is to prove $W_n \vertk W$ for some centred
Gaussian random element of ${\cal H}$. By the Continuous Mapping Theorem (CMT), we would then have $T_n \vertk \|W\|^2_\oplus$. Here and in the sequel,
$\vertk$ denotes convergence in distribution of random elements (especially: of random variables). Moreover,
$o_\PP(1)$ refers to convergence in probability to zero of random elements ${\cal H}$.
The main result of this section is as follows.

\begin{thm} {\rm{(}}Convergence of $W_n${\rm{)}} \label{thmmain}\\
Suppose that $X$ has some non-degenerate $d$-variate normal distribution. If the weight function $w_\gamma$ in (\ref{weightf}) satisfies $\gamma >2$,
 there is some centred Gaussian random element $W$ of ${\cal H}$ having
covariance {\rm{(}}matrix{\rm{)}} kernel $K(s,t) = \BE\big{[} W(s) W(t)^\top \big{]}$, where
\begin{equation}\label{kernel}
K(s,t) = {\rm{e}}^{(\|s\|^2+ \|t\|^2)/2} \left({\rm{e}}^{s^\top t} \left(ts^\top + \rm{I}_d\right)
- ts^\top - (1+  s^\top t) \, {\rm{I}}_d \right),
\end{equation}
$s,t \in \RR^d$, so that $W_n \vertk W$ as $n \to \infty$.
\end{thm}

\vspace*{3mm}

\begin{cor} \label{cormain}
The limit distribution of $T_{n,\gamma}$ as $n \to \infty$ under $H_0$ is that of
\[
 \|W\|^2_\oplus = \int_{\RR^d} \|W(t)\|^2 \, w_\gamma(t)   \, {\rm{d}}t.
\]
\end{cor}

\vspace*{3mm}

\noindent {\sc Proof} of Theorem \ref{thmmain}.
To highlight the main idea of the proof, let
\[
Z_j(t) := \e^{t^\top Y_{n,j}}\left(Y_{n,j} - t\right), \qquad Z_j^\circ(t) := \e^{t^\top X_j}\left(X_j - t\right), \qquad t \in \RR^d,
\]
and put
\begin{equation}\label{wnull}
W_n^\circ(t) :=  \frac{1}{\sqrt{n}} \sum_{j=1}^n Z_j^\circ(t).
\end{equation}
Notice that $Z_1^\circ,Z_2^\circ, \ldots $ are i.i.d. (centred) random elements of ${\cal H}$. Moreover, the condition $\gamma >2$ implies  $\BE \|Z_1^\circ\|^2_\oplus < \infty$.
 By a Hilbert space central limit theorem (CLT), see e.g. \cite{KMM2000}, there is some centred Gaussian random element $W^\circ$ of ${\cal H}$ such that
$W_n^\circ \vertk W^\circ$. The idea now is to approximate $W_n$ by a suitable random element $\widetilde{W}_n$ of ${\cal H}$ so that
$\|W_n - \widetilde{W}_n\|_\oplus = o_\PP(1)$ as $n \to \infty$, and
$\widetilde{W}_n(t) = n^{-1/2} \sum_{j=1}^n \widetilde{Z}_j(t)$,
where $\widetilde{Z}_1,\widetilde{Z}_2, \ldots $ are i.i.d. centred random elements of ${\cal H}$ satisfying $\BE \|\widetilde{Z}_1\|^2_\oplus < \infty$.
The assertion would then follow from the CLT in Hilbert spaces and Slutzky's lemma.
To this end, put
\begin{equation}\label{delta1}
\Delta_{n,j} = Y_{n,j} - X_j = \left(S_n^{-1/2}- \textrm{I}_d\right)X_j - S_n^{-1/2}\overline{X}_n, \qquad j=1,\ldots,n.
\end{equation}
A Taylor expansion gives
\begin{equation}\label{Taylor}
\e^{t^\top \Delta_{n,j}}  =  1 + t^\top \Delta_{n,j} + \frac{1}{2} \left(t^\top \Delta_{n,j}\right)^2 \exp\left(\Theta_{n,j} t^\top \Delta_{n,j} \right),
\end{equation}
where $|\Theta_{n,j}| \le 1$. By some algebra, it follows that
\begin{equation}\label{zerleg}
W_n(t)  =  \frac{1}{\sqrt{n}}\sum_{j=1}^n \e^{t^\top X_j} (X_j-t)  + A_n(t) + B_n(t) + C_n(t),
\end{equation}
where
\begin{eqnarray} \label{ant}
A_n(t) & = & \frac{1}{\sqrt{n}}\sum_{j=1}^n \e^{t^\top X_j} t^\top \Delta_{n,j} \,(X_j-t),\\ \label{bnt}
B_n(t) & = &  \frac{1}{\sqrt{n}}\sum_{j=1}^n \e^{t^\top X_j} \left(1+ t^\top \Delta_{n,j}\right) \Delta_{n,j},\\ \label{cnt}
C_n(t) & = & \frac{1}{\sqrt{n}}\sum_{j=1}^n \e^{t^\top X_j} \, \frac{1}{2} \, \left(t^\top \Delta_{n,j}\right)^2 \exp\left(\Theta_{n,j} t^\top \Delta_{n,j}\right) \left(X_j-t+\Delta_{n,j}\right).
\end{eqnarray}
Notice that the first term on the right-hand side of (\ref{zerleg}) is $W_n^\circ(t)$, as  given in (\ref{wnull}).
By Proposition \ref{propcnt}, we have $\|C_n\|_\oplus = o_\PP(1)$, and Proposition \ref{propbnt} yields
\[
B_n(t) = - \frac{1}{\sqrt{n}} \, \exp \left(\frac{\|t\|^2}{2}\right) \sum_{j=1}^n \left( \frac{X_jX_j^\top t - t}{2} + X_j \right) + o_\PP(1).
\]
In view of  Proposition \ref{propant},  (\ref{zerleg})
implies
$W_n(t) = n^{-1/2} \sum_{j=1}^n \widetilde{Z}_j(t) + o_\PP(1)$,
where
\begin{equation} \label{widetildez}
\widetilde{Z}_j(t) = \e^{t^\top X_j} (X_j-t) - \e^{\|t\|^2/2}  \left( X_jX_j^\top t - t + X_j \right).
\end{equation}
A straightforward calculation shows $\BE \widetilde{Z}_j(t) = 0$, $t \in \RR^d$. Moreover,  due to the condition $\gamma >2$, we
have $\BE \|\widetilde{Z}_j\|^2_\oplus < \infty$. Hence, $\widetilde{Z}_1, \widetilde{Z}_2, \ldots $ are i.i.d. centred random elements of ${\cal H}$,
and the CLT in Hilbert spaces yields the assertion, since the kernel $K$ figuring in (\ref{kernel}) is given by $K(s,t) = \BE[ \widetilde{Z}_1(s) \widetilde{Z}_1(t)^\top ]$
(for details in computing $K(s,t)$, see Proposition \ref{calckernel}). \bewend

\vspace*{3mm}

The following result provides some information on the limit null distribution of $T_{n,\gamma}$.

\begin{thm} Let $T_{\infty,\gamma}$ be a random variable with the distribution of $\|W\|^2_\oplus$,  where $W$ is given in Theorem \ref{thmmain}.
We then have
\[
\BE\left[ T_{\infty,\gamma}\right] = \left( \frac{\pi}{\gamma - 2}\right)^{d/2} \left( d + \frac{d}{2(\gamma -2)}\right) - \left(\frac{d(d+1)}{2(\gamma -1)} + d\right)\,
\left( \frac{\pi}{\gamma - 1}\right)^{d/2} .
\]
\end{thm}

\noindent {\sc Proof.} By Fubini's theorem, it follows that
\[
\BE\left[ T_{\infty,\gamma}\right] =  \int_{\RR^d} \BE \|W(t)\|^2 \, w_\gamma(t) \, \textrm{d}t.
\]
Writing $\textrm{tr}(D)$ for the trace of a square matrix $D$, we have
\begin{eqnarray*}
\BE \|W(t)\|^2  & = & \BE \left[W(t)^\top W(t) \right] \ = \ \BE \left[ \textrm{tr}(W(t)^\top W(t))\right] \\
& = & \BE\left[\textrm{tr}(W(t) W(t)^\top)\right] \ = \ \textrm{tr} \left(\BE\left[W(t) W(t)^\top \right] \right) \ = \
 \textrm{tr} \left(K(t,t)\right)\\
& = & \e^{\|t\|^2} \left( \e^{\|t\|^2} \left(d+ \|t\|^2\right) - \left( \|t\|^2 + d \|t\|^2\right) - d \right).
\end{eqnarray*}
Now, some straightforward algebra yields the assertion. \bewend

In the univariate case $d=1$, we also obtained an explicit expression for the variance of $T_{\infty,\gamma}$, by using the relation
\[
\textrm{Var}[T_{\infty,\gamma}] =  2 \iint_{\mathbb{R}^2} K^2(s,t) w_\gamma(s) w_\gamma(t) \, \textrm{d}s \textrm{d}t.
\]
By tedious calculations, it follows that
\[
\textrm{Var}[T_{\infty,\gamma}]  = 2\pi \left( \beta^{-1} + \beta^{-3} + \delta + \delta^3 + \frac{1}{4} (\beta^2 +2)\delta^5
 -4\eta -12 \eta^3 -16(2\beta^2+1) \eta^5 \right),
\]
where $\beta = \gamma -1$, $\delta = (\beta^2-1)^{-1/2}$ and $\eta = (4\beta^2 -1)^{-1/2}$.

 Table 1 contains expectation and variance of $T_{\infty,\gamma}$ in the univariate case for various values of $\gamma$.

%% For tables use
%\begin{center}
%\begin{table}
%\begin{center}
%% table caption is above the table
%\caption{Expectation and variance of the limit null distribution when $d=1$}
%\label{tab0}       % Give a unique label
%% For LaTeX tables use
%\begin{tabular}{|c|cc|}
%\hline
%\hline
%$\gamma$ & $\mathbb{E}[T_{\infty,\gamma}]$ & $\textrm{Var}[T_{\infty,\gamma}]$\\ % & $\textrm{Skew}[T_{\infty,\gamma}]$ \\
%\hline
%$2.5$ & $2.6013$ & $4.7153$\\% & $1.9576$ \\
%$3$ &    $0.7787$ & $0.5430$\\% & $2.1780$ \\
%$4$ &    $0.2022$ & $0.0458$\\% & $2.4423$ \\
%$5$ &    $0.0861$ & $0.0093$\\% & $2.5812$ \\
%$7$ &    $0.0277$ & $0.0011$\\% & $2.7053$ \\
%%$10$ &  $0.0094$ & $0.0001$ & $2.7700$ \\
%\hline
%\end{tabular}
%\end{center}
%\end{table}
%\end{center}

\begin{center}
\begin{table}
\begin{center}
%\label{tab0}       % Give a unique label
% For LaTeX tables use
\begin{tabular}{|c|ccccc|}
\hline
\hline
$\gamma$ &  2.5 & 3 & 4 & 5 & 7\\
\hline
$\mathbb{E}[T_{\infty,\gamma}]$ & 2.6013 & 0.7787 & 0.2022 & 0.0861 & 0.0277 \\
$\textrm{Var}[T_{\infty,\gamma}]$ & 4.7153 & 0.5430 & 0.0458 & 0.0094 & 0.0011\\ % & $\textrm{Skew}[T_{\infty,\gamma}]$ \\
\hline
\end{tabular}
\end{center}
\begin{center}
{\bf Table 1:} Expectation and variance of the limit null distribution when $d=1$
\end{center}
\end{table}
\end{center}

%%%%%%%%%%%%%%%%%%%%%%%%%%%%%%%%%%%%%%%%%%%%%%%%%%%%%%%%%%%%%%%%%%%%%%%%%%%%%%%%%%%%%%%%%%%%%%%%%%%%%%%%%%%%%%%%%%%%%%%%%%%%%%%%%%%%%%%%%%%%%%%%%%%%%%%%%%%%%%%%%%%%%%%%%%%%%%
%
%
%                   Consistency
%
%
%%%%%%%%%%%%%%%%%%%%%%%%%%%%%%%%%%%%%%%%%%%%%%%%%%%%%%%%%%%%%%%%%%%%%%%%%%%%%%%%%%%%%%%%%%%%%%%%%%%%%%%%%%%%%%%%%%%%%%%%%%%%%%%%%%%%%%%%%%%%%%%%%%%%%%%%%%%%%%%%%%%%%%%%%%%%%%

\section{Consistency}\label{secconsist}
In this section, let $X$ have an absolutely continuous distribution, and suppose that $m(t) = \BE [\exp(t^\top X)]$ exists for each $t \in \RR^d$.
In view of affine invariance, let w.l.o.g. $\BE[X] = 0$ and $\BE [XX^\top] = \textrm{I}_d$.
\begin{thm}\label{thmconsist}
We have
\begin{equation}\label{righthside}
\liminf_{n\to \infty} \frac{T_{n,\gamma}}{n} \ge \int_{\RR^d} \|m'(t) - tm(t)\|^2 \, w_\gamma(t) \, {\rm{d}}t \quad \PP\textrm{-almost surely}.
\end{equation}
\end{thm}

\vspace*{3mm}

\noindent {\sc Proof}. Fix $K>0$, and put $S(K) := \{t \in \RR^d:\|t\| \le K\}$. Let $B(K)$ be the Banach space of continuous functions
$f:S(K) \rightarrow \RR$, equipped with the norm $\|f\|_\infty = \max_{\|t \|\le K} |f(t)|$. Recall $M_n(t)$ from (\ref{defmgf}),
and put
\[
M_{n,0}(t) = \frac{1}{n}\sum_{j=1}^n \textrm{e}^{t^\top X_j}, \quad t \in \RR^d.
\]
Let $\xi_n = \max_{j=1,\ldots,n} \|\Delta_{n,j}\|$, where $\Delta_{n,j}$ is given in (\ref{delta1}). From (\ref{delta1}) we obtain
\begin{equation}\label{xinungl}
\xi_n \le \big{\|} S_n^{-1/2} - \textrm{I}_d\|_2 \cdot \max_{1\le j \le n}\|X_j\| + \| S_n^{-1/2}\|_2 \cdot \|\overline{X}_n\|.
\end{equation}
Since the existence of $m(t)$ implies $\BE\|X\|^4 < \infty$, Theorem 5.2 of \cite{BAR1963} yields
\begin{equation}\label{max14}
n^{-1/4} \max_{1\le j \le n}\|X_j\| \to 0 \quad \PP\textrm{-almost surely.}
\end{equation}
From $S_n - \textrm{I}_d = n^{-1}\sum_{j=1}^n (X_jX_j^\top- \textrm{I}_d) - \overline{X}_n\overline{X}_n^\top$ and
Kolmogorov's variance criterion for averages (see \cite{KAL2002}, p. 73), we have
$\lim_{n\to \infty} n^{1/2-\varepsilon}\|S_n- \textrm{I}_d\|_2 = 0$ $\PP\textrm{-a.s.}$
for each $\varepsilon >0$. Hence Proposition \ref{proplinalgebra} implies
\begin{equation}\label{limxin}
\lim_{n \to \infty} \xi_n = 0 \quad \PP\textrm{-a.s.}
\end{equation}
Using the notation $\|f\|_\infty = \max_{\|t \|\le K}\|f(t)\|$ also for a function $f:\RR^d \to \RR^d$, (\ref{Taylor}) gives
\begin{equation}\label{mnungl}
\|M_n - M_{n,0}\|_\infty \le \|M_{n,0}\|_\infty \cdot \xi_n \cdot K \cdot \left(1+ \frac{K\xi_n}{2}\textrm{e}^{K\xi_n}\right).
\end{equation}
By the strong law of large numbers (SLLN) in $B(K)$, the first factor on the right-hand side converges almost surely to $\|m\|_\infty$,
and thus (\ref{limxin}) entails
\begin{equation}\label{ascon1}
\lim_{n \to \infty} \|M_n - M_{n,0}\|_\infty = 0 \quad \PP\textrm{-a.s.}
\end{equation}
Putting
\begin{equation}\label{defmaxnorm}
F_n = \max_{j=1,\ldots,n} \|X_j\|,
\end{equation}
the triangle inequality gives
\[
\|M_n'-M_{n,0}'\|_\infty \le \left(\xi_n + F_n \right) \cdot \|M_n-M_{n,0}\|_\infty + \xi_n \cdot \|M_{n,0}\|_\infty.
\]
Invoking (\ref{mnungl}), (\ref{xinungl}), (\ref{max14}) and Proposition \ref{proplinalgebra}, we have
\begin{equation}\label{ascon2}
 \lim_{n\to \infty} \|M_n'-M_{n,0}'\|_\infty = 0 \quad \PP\textrm{-a.s.}
\end{equation}
Writing id for the identity function on $\RR^d$,  the triangle inequality yields
\[
\|M_n^{'}  - \textrm{id}\cdot M_n\|_\infty \le \|M_n'-M_{n,0}'\|_\infty + \|M_{n,0}'-\textrm{id}\cdot M_{n,0}\|_\infty + K\cdot \|M_{n,0}-M_n\|_\infty.
\]
From (\ref{ascon1}), (\ref{ascon2}) and the SLLN in $B(K)$  it follows that
\begin{equation}\label{ungl1}
\limsup_{n\to \infty} \|M_n^{'}  - \textrm{id}\cdot M_n\|_\infty \le \|m' - \textrm{id}\cdot m\|_\infty  \quad \PP\textrm{-a.s.}
\end{equation}
Likewise, we have
\[
\|M_{n,0}^{'}  - \textrm{id}\cdot M_{n,0}\|_\infty \le \|M_{n,0}'-M_n'\|_\infty + \|M_n'-\textrm{id}\cdot M_n\|_\infty + K\cdot \|M_n-M_{n,0}\|_\infty
\]
and thus
\begin{equation}\label{ungl2}
\|m' - \textrm{id}\cdot m\|_\infty \le \liminf_{n\to \infty} \|M_n^{'}  - \textrm{id}\cdot M_n\|_\infty \quad \PP\textrm{-a.s.}
\end{equation}
Upon combining (\ref{ungl1}) and (\ref{ungl2}), and using Fatou's lemma, we obtain
\begin{eqnarray*}
\liminf_{n\to \infty} \frac{T_{n,\gamma}}{n} & \ge & \liminf_{n \to \infty} \int_{S(K)} \|M_n'(t) - tM_n(t)\|^2 \, w_\gamma(t)\, \textrm{d}t \\
& \ge & \int_{S(K)} \liminf_{n\to \infty} \|M_n'(t) - tM_n(t)\|^2 \, w_\gamma(t)\, \textrm{d}t \\
& = & \int_{S(K)} \|m'(t) - t m(t)\|^2 \, w_\gamma(t)\, \textrm{d}t \quad \PP\textrm{-a.s.}
\end{eqnarray*}
Since $K$ was arbitrary, the assertion follows. \bewend

\vspace*{3mm}

\begin{rem}
If $X$ has some non-degenerate non-normal distribution with existing moment generating function,
then $m'(t) \neq tm(t)$ for at least one $t$. Consequently, the right-hand side of (\ref{righthside}) is strictly positive, and thus
\begin{equation}\label{conjecturinf}
\lim_{n\to \infty} T_{n,\gamma} = \infty \quad \PP\textrm{-a.s.}
\end{equation}
Hence, the test is consistent against each such alternative. We conjecture that (\ref{conjecturinf}) holds for {\em any}
non-normal alternative distribution. A proof of such a result, however, remains an open problem.
\end{rem}

\vspace*{3mm}

%%%%%%%%%%%%%%%%%%%%%%%%%%%%%%%%%%%%%%%%%%%%%%%%%%%%%%%%%%%%%%%%%%%%%%%%%%%%%%%%%%%%%%%%%%%%%%%%%%%%%%%%%%%%%%%%%%%%%%%%%%%%%%%%%%%%%%%%%%%%%%%%%%%%%%%%%%%%%%%%%%%%%%%%%%%%%%
%
%
%                   The limit gamma to infinity
%
%
%%%%%%%%%%%%%%%%%%%%%%%%%%%%%%%%%%%%%%%%%%%%%%%%%%%%%%%%%%%%%%%%%%%%%%%%%%%%%%%%%%%%%%%%%%%%%%%%%%%%%%%%%%%%%%%%%%%%%%%%%%%%%%%%%%%%%%%%%%%%%%%%%%%%%%%%%%%%%%%%%%%%%%%%%%%%%%

\section{The limit $\gamma \to \infty$}\label{seclimgam}
In this section, we will show that, for fixed $n$, the statistic $T_{n,\gamma}$, after a suitable scaling, approaches a linear combination of two well-known
measures of multivariate skewness as $\gamma \to \infty$.

\begin{thm}\label{gammainf}
We have {\rm{(}}pointwise on the underlying probability space{\rm{)}}
\[
\lim_{\gamma \to \infty} \gamma^{2+d/2} \, \frac{16 T_{n,\gamma}}{n \pi^{d/2}} = 2b_{1,d} +  \widetilde{b}_{1,d},
\]
where
\begin{equation}\label{a.b}
b_{1,d} = \frac{1}{n^2}  \sum_{j,k=1}^n (Y_{n,j}^\top Y_{n,k})^3
\end{equation}
is nonnegative invariant sample skewness in the sense of \cite{Mar1970},
and
\[
\widetilde{b}_{1,d} =  \frac{1}{n^2} \sum_{j,k=1}^n Y_{n,j}^\top Y_{n,k} \, \|Y_{n,j}\|^2 \, \|Y_{n,k}\|^2
\]
denotes sample skewness in the sense of M\'{o}ri, Rohatgi and Sz\'{e}kely, as defined in \cite{MRS1993}.
\end{thm}

\noindent {\sc Proof.} We start with (\ref{rep1}) and use
\[
\exp\left(\frac{\|Y_{n,j}+Y_{n,k}\|^2}{4 \gamma}\right) = 1 + \frac{\|Y_{n,j}+Y_{n,k}\|^2}{4\gamma} + \frac{\|Y_{n,j}+Y_{n,k}\|^4}{32\gamma^2} + O\left(\gamma^{-3}\right)
\]
as $\gamma \to \infty$. Multiplying this expression with the term within the rightmost bracket of (\ref{rep1}), and using the relations $\sum_{j=1}^n Y_{n,j}=0$, $\sum_{j=1}^n \|Y_{n,j}\|^2=nd$,
$\sum_{j,k=1}^n \|Y_{n,j}+ Y_{n,k}\|^2  =  2n^2d$,
\begin{eqnarray*}
\sum_{j,k=1}^n \|Y_{n,j}+ Y_{n,k}\|^4 & = & 2n^2 \left(\frac{1}{n} \sum_{j=1}^n \|Y_{n,j}\|^4 + d^2 + 2d\right), \\
\sum_{j,k=1}^n \|Y_{n,j}+ Y_{n,k}\|^4 Y_{n,j}^\top Y_{n,k} & = & 8 \sum_{j,k=1}^n \left(Y_{n,j}^\top Y_{n,k}\right)^2 \|Y_{n,j}\|^2 + 4n^2 b_{1,d} + 2 n^2 \widetilde{b}_{1,d}
\end{eqnarray*}
as well as
\begin{eqnarray*}
\sum_{j,k=1}^n \left(Y_{n,j}^\top Y_{n,k}\right)^2 \|Y_{n,j}\|^2  & = & \sum_{j,k=1}^n \textrm{tr}\left(Y_{n,j}^\top Y_{n,k}Y_{n,k}^\top Y_{n,j} \, \|Y_{n,j}\|^2\right) \\
& = & \textrm{tr}\left( \sum_{k=1}^n Y_{n,k}Y_{n,k}^\top \sum_{j=1}^n Y_{n,j}Y_{n,j}^\top \|Y_{n,j}\|^2\right)\\
& = & \textrm{tr}\left( n \textrm{I}_d \sum_{j=1}^n Y_{n,j}Y_{n,j}^\top \|Y_{n,j}\|^2\right)\\
& = & n \, \sum_{k=1}^n \textrm{tr}\left(Y_{n,j}^\top Y_{n,j} \|Y_{n,j}\|^2\right)\\
& = & n \sum_{j=1}^n \|Y_{n,j}\|^4,
\end{eqnarray*}
the result follows by tedious but straightforward algebra.
\bewend

\begin{rem}
It is interesting to note a similarity between Theorem \ref{gammainf} and Theorem 2.1 of \cite{Hen1997}, who showed that the
BHEP-statistic for testing for multivariate normality, after suitable rescaling, approaches the linear combination $2b_{1,d} + 3 \widetilde{b}_{1,d}$ as a smoothing parameter
(called $\beta$ in that paper) tends to $0$. Since $\beta$ and $\gamma$ are related by $\beta = \gamma^{-1/2}$, this corresponds to letting $\gamma$ tend to infinity.
The same linear combination $2b_{1,d} + 3 \widetilde{b}_{1,d}$ also showed up as a limit statistic in \cite{HJM2017}. Notice that, in the univariate case, the
limit statistic figuring in Theorem \ref{gammainf} is nothing but three times squared sample skewness.
It should be stressed that tests for multivariate normality based on $b_{1,d}$ or $\widetilde{b}_{1,d}$ {\rm{(}}or on related measures of multivariate skewness and kurtosis{\rm{)}}
lack consistency against general alternatives, see, e.g., \cite{BH1991}, \cite{BH1992},  \cite{Hen1994a}, and \cite{Hen1994b}.
\end{rem}

%%%%%%%%%%%%%%%%%%%%%%%%%%%%%%%%%%%%%%%%%%%%%%%%%%%%%%%%%%%%%%%%%%%%%%%%%%%%%%%%%%%%%%%%%%%%%%%%%%%%%%%%%%%%%%%%%%%%%%%%%%%%%%%%%%%%%%%%%%%%%%%%%%%%%%%%%%%%%%%%%%%%%%%%%%%%%%
%
%
%                   Monte Carlo results
%
%
%%%%%%%%%%%%%%%%%%%%%%%%%%%%%%%%%%%%%%%%%%%%%%%%%%%%%%%%%%%%%%%%%%%%%%%%%%%%%%%%%%%%%%%%%%%%%%%%%%%%%%%%%%%%%%%%%%%%%%%%%%%%%%%%%%%%%%%%%%%%%%%%%%%%%%%%%%%%%%%%%%%%%%%%%%%%%%

\section{Monte Carlo results}\label{secsimul}

In this section we compare the finite-sample power performance of the newly proposed test to those of several competing tests for normality,
both for the univariate and the multivariate case.
In the case $d=1$ the competing procedures are
\begin{enumerate}
%\item[a)] the Lilliefors ($Lil$), also referred to as the Kolmogorov-Smirnov, test,
\item[a)] the Cram\'{e}r-von Mises ($CvM$) test,
\item[b)] the Anderson-Darling ($AD$) test,
\item[c)] the Shapiro-Wilk ($SW$) test,
\item[d)] the Jarque-Bera ($JB$) test,
\item[e)] the Zghoul ($Z$) test.
\end{enumerate}

The R package $nortest$ contains the functions cvm.test and ad.test, which can be used in order to calculate the test statistic
 and the associated p-value for each of the first two tests mentioned above, see \cite{GL2015}. The Shapiro-Wilk test can be performed
 using the $Shapiro.test$ function included in the $stats$ package. The R package $tseries$ contains a function called $jarque.bera.test$, which
 can be used in order to calculate the test statistic and p-value associated with the fourth test mentioned above, see \cite{TH2017}.
Each of these tests are well-known and will not be discussed further.

The test of Zghoul (see \cite{Zgh2010}) is based on the empirical moment generating function. \cite{Zgh2010} includes a Monte Carlo study which indicates that the finite-sample power performance of the test compares favorably to that of its competitors, especially against symmetric alternatives with kurtosis higher than that of the normal distribution. However, the mentioned paper fails to provide the mathematical theory underlying this test. \cite{HEK2017} more recently provided this theory, including a proof that the test is consistent against general alternatives.

The test statistic of Zghoul is a weighted $L^2$-statistic, given by
\[
Z_{n} (\gamma) = n \int_{\mathbb{R}} (M_n(t) - m(t))^2 \textrm{exp}(-\gamma t^2) \textrm{d}t,
\]
where $\gamma>2$ is a smoothing parameter. Based on the finite-sample performance reported in \cite{Zgh2010}, the author recommended setting $\gamma$ equal to either $3$ or $15$. The numerical results presented below include the powers obtained using both of these values for the smoothing parameter; the resulting tests are denoted by $Z_{3}$ and $Z_{15}$ respectively. The test statistic $Z_{n}(\gamma)$
can be rewritten in the computationally amenable form
\[
Z_{n} (\gamma) = \sqrt \pi \left[ \frac{n}{\sqrt{\gamma \! - \! 1}} - \frac{2}{\sqrt{\gamma \! -\!  \frac{1}{2}}}\! \sum_{i=1 }^n
\textrm{exp} \! \left( \! \frac{Y^2 _{n,i}}{4 \gamma - 2}\right)
 + \frac{1}{n \sqrt{\gamma}} \sum_{i,j=1}^n \textrm{exp}  \! \left(\!  \frac{\left( {Y _{n,i}} + {Y _{n,i}} \right)^2}{4 \gamma} \right) \right].
\]
This test rejects normality for large values of the test statistic.

Turning our attention to the multivariate case, we consider the finite-sample power performance of the newly proposed test to
those of some prominent competing tests, and to two very recent tests for multinormality.
These procedures are
\begin{enumerate}
\item[a)] Mardia's tests based on skewness and kurtosis,
\item[b)] the Henze-Zirkler test,
\item[c)] the {\em energy test} of Sz\'{e}kely and Rizzo,
\item[d)] a recent test of Henze, Jim\'{e}nez-Gamero and Meintanis,
\item[e)] a recent test of Henze and Jim\'{e}nez-Gamero.
\end{enumerate}

\vspace*{1mm}

\noindent {\bf a):} {\bf Mardia's tests based on skewness and kurtosis}\\[1mm]
Mardia's test for multinormality based on sample skewness rejects $H_0$ for large values of $b_{1,d}$, where $b_{1,d}$ is given in
(\ref{a.b}).  Notice that $b_{1,d}$ is a consistent estimator of
$\beta_{1,d} = E (X_1^\top X_2)^3$. Under normality we have $\beta_{1,d}=0$,
 and the limit distribution of $nb_{1,d}$ as $n \rightarrow \infty$ is $6 \chi^2_{d(d+1)(d+2)/6}$, see \cite{Mar1970}. The limit distribution
 of $nb_{1,d}$ under elliptical symmetriy has been derived by \cite{BH1992}.

Sample kurtosis in the sense of Mardia is given by
\[
b_{2,d} = \frac{1}{n}  \sum_{j=1}^n \|Y_{n,j}\|^4,
\]
which is an estimator of $\beta_{2,d} = \mathbb{E}\|X_{1}\|^4$.

Under normality, we have $\beta_{2,d}=d(d+2)$, and the limit null distribution of $\sqrt{n} (b_{2,d}-d(d+2))$ is the normal distribution $\textrm{N}(0,8d(d+2))$, see \cite{Mar1970}. The test based on
$b_{2,d}$ rejects $H_0$ for large and small values of $\beta_{2,d}$.

The R package $QuantPsyc$ contains a function $mult.norm$, which calculates both Mardia's skewness and kurtosis measures as well as the p-values associated with the corresponding tests from multivariate normality, see \cite{Fle2012}. Below we denote the tests based on Mardia's skewness and kurtosis by $MS_n$ and $MK_n$, respectively.

We stress that there are several other measures of skewness and kurtosis, see Sections 3 and 4 of \cite{Hen2002}. The deficiences of such measures
as statistics for purportly ``directed  tests'' for multivariate normality have been addressed in \cite{BH1991}, \cite{BH1992}, \cite{Hen1994a} as well as \cite{Hen1994b}.

\vspace*{2mm}

\noindent {\bf b):} {\bf The Henze-Zirkler test}\\[1mm]
Writing $\Psi_n(t) = n^{-1}\sum_{k=1}^n \exp(\textrm{i}t^\top Y_{n,k})$ for the empirical {\em characteristic} function of the scaled residuals $Y_{n,1},\ldots,Y_{n,n}$,
\cite{HZ1990} proposed the test statistic
\[
HZ_{n} (\gamma) = (2 \pi \gamma^2)^{-d/2} \int_{\mathbb{R}^d} \bigg{|} \Psi_n(t) - \exp \left( - \frac{\|t|^2}{2} \right) \bigg{|}^2 \,
\exp \left( -\frac{\|t\|^2}{2 \gamma^2} \right)  \, \textrm{d}t,
\]
for some fixed $\gamma>0$. The test statistic can be rewritten as
\begin{eqnarray*}
HZ_{n}(\gamma) &=& \frac{1}{n^2}  \sum_{j,k=1}^n \exp\left(-\frac{\gamma^2}{2} \|Y_{n,j}-Y_{n,k}\|^2\right) \\
& & \quad \quad - 2(1+\gamma^2)^{-d/2} \frac{1}{n} \sum_{j=1}^n \exp \left( -\frac{\gamma^{2} \|Y_{n,j}\|^{2}}{2(1+\gamma^2)} \right)
+ (1+2\gamma^2)^{-d/2}.
\end{eqnarray*}
The Henze-Zirkler test is obtained by setting
\[
\gamma = \frac{1}{\sqrt{2}}  \left( \frac{(2d+1)n}{4} \right)^{1/(d+4)}.
\]
This choice of $\gamma$ corresponds to the optimal bandwidth for a multivariate nonparametric density estimator with a Gaussian kernel.
The Henze-Zirkler test is included because of its impressive power performance reported in previous studies, see, e.g.,  \cite{MM2005}.

The Henze-Zirkler test (denoted $HZ_n$ below) is programmed in the function $hzTest$ in the R package $MVN$, see \cite{KGZ2014}.
This test rejects normality for large values of the test statistic. $HZ_n$ is equivalent to a test by \cite{BF1993}, as remarked in \cite{Hen2002}.

\vspace*{2mm}

\noindent {\bf c):} {\bf The energy test}\\[1mm]
\cite{SR2005} proposed the test statistic
\[
EN_{n} = n \left(\! \frac{2}{n}  \sum_{j=1}^n \mathbb{E} \|Y_{n,j}\! -\! Z\| - \frac{2 \Gamma((d+1)/2)}{\Gamma(d/2)}
- \frac{1}{n^2} \sum_{j,k=1}^n \mathbb{E} \|Y_{n,j}\! -\! Y_{n,k}\| \! \right).
\]
Here, the first expectation is with respect to a random vector $Z$, which is independent of $Y_{n,j}$ and has the distribution $\textrm{N}_d(0,\textrm{I}_d)$, and
\begin{eqnarray*}
\mathbb{E} \|a-Z\|=\sqrt{2} \frac{2 \Gamma((d+1)/2)}{\Gamma(d/2)} &+& \sqrt\frac{2}{\pi}\sum_{k=0}^\infty \left\{ \frac{(-1)^k}{k!2^k}
\frac{\|a\|^{2k+2}}{(2k\! +\! 1)(2k\! +\! 2)} \right.
\\
& & \left. \quad \quad \times \frac{2 \Gamma((d+1)/2)\Gamma(k\! +\! 1.5)}{\Gamma((d/2)\! +\! k\! +1)} \right\}.
\end{eqnarray*}
This test, denoted by $EN_n$, is known as the energy test. Rejection of $H_0$ is for large values of $EN_n$. The function
 $mvnorm.etest$  in the R package $energy$ calculates this test statistic as well as the corresponding p-value, see \cite{RS2016}.
The energy test is also reported to have excellent power performance, see \cite{JV2014}.

\vspace*{2mm}

\noindent {\bf d):} {\bf The Henze--Jim\'{e}nez-Gamero--Meintanis test}\\[1mm]
By generalizing a characterization of normality involving both the characteristic and the moment generating function (see \cite{Vol2014}) to the multivariate
case, \cite{HJM2017} proposed to base a test of $H_0$ on the weighted $L^2$-statistic
\begin{eqnarray*}
HM_n=n \int_{\mathbb{R}^d} \! \left(\!  \frac{1}{n^2} \sum_{j=1}^n \cos \left( t^\top Y_{n,j}\right) \sum_{j=1}^n \exp \left( t^\top Y_{n,j}\right) -1\! \right)^2 \exp \left( -\gamma \|t\|^2 \right) \textrm{d}t,
\end{eqnarray*}
for some parameter $\gamma>1$. Putting $Y_{jk}^\pm=Y_{n,j} \pm Y_{n,k}$, $HM_n$ can be rewritten as
\begin{eqnarray*}
HM_n &=& \left( \frac{\pi}{\gamma} \right)^{d/2} \! \left\{ { \! \frac{1}{2 n^3} \sum_{j,k,\ell,m=1}^n \left\{ {\exp \left(
\frac{\| Y_{jk}^+ \|^2 - \| Y_{\ell m}^- \|^2}{4 \gamma} \right) \cos\left( \frac{Y_{jk}^{+ \top} Y_{ \ell m}^-}{2 \gamma} \right) }\right. }\right.
\\
& & \left. { + \exp \left( \frac{\| Y_{jk}^+ \|^2 - \| Y_{\ell m}^+ \|^2}{4 \gamma} \right) \cos \left( \frac{Y_{jk}^{+ \top} Y_{\ell m}^+}{2 \gamma} \right)
}\right\} \\
& & \left. { - \frac{2}{n} \sum_{j,k=1}^n \exp \left( \frac{\| Y_{n,j} \|^2 - \| Y_{n,k} \|^2}{4 \gamma}
\cos \left( \frac{Y_{n,j}^\top Y_{n,k}}{2 \gamma} \right) +n \right) }\right\}.
\end{eqnarray*}
The test based on $HM_n$ uses an upper rejection region.

A disadvantage of this test is the need to calculate a four-fold sum in order to evaluate the test statistic.
 Hence, the amount of computer time required in order to perform this test is substantially more than the time required for the other tests under consideration.

\vspace*{2mm}

\noindent {\bf e):} {\bf The Henze--Jim\'{e}nez-Gamero test}\\[1mm]
\cite{HJ2017} present a multivariate generalization of a recent class of tests for univariate normality (see \cite {HEK2017}) based on the empirical moment generating function. The test statistic is
\begin{eqnarray*}
HJ_n=n \int_{\mathbb{R}^d} \left( M_n (t)-m(t)\right) ^2 \exp \left (-\beta \|t\|^2 \right) \textrm{d}t,
\end{eqnarray*}
where $\beta>1$ is a fixed parameter. An alternative representation of $HJ_n$ is
\begin{eqnarray*}
HJ_n &=& \pi ^{d/2} \left\{ \frac{1}{n} \sum_{j,k=1}^n \frac{1}{\beta^{d/2}} \exp \left( \frac{\| Y_{n,j}+Y_{n,j} \|^2}{4\beta} \right) + \frac{n}{(\beta - 1)^{d/2}} \right.\\
& & \left. { -2 \sum_{j=1}^n \frac{1}{(\beta-1/2)^{d/2}} \exp \left( \frac{\|Y_{n,j}\|^2}{4\beta-2} \right) }\right\},
\end{eqnarray*}
which is amenable to computation. This test rejects $H_0$ for large values of the test statistic.

The numerical results below were obtained using the software package R, see \cite{R2015}.
Monte Carlo simulation was used in order to estimate the upper percentiles of the distribution of $T_{n,\gamma}$ for the sample sizes $n=20,50,100,200,300$
and dimensions $d=1,2,3,5$. Table 2 provides the empirical $95$-percentiles of $(16 \gamma^{2+d/2}/ \pi^{d/2})T_{n,\gamma}$ obtained
by simulating one million samples from a $\textrm{N}_d(0,\textrm{I}_d)$ distribution using $\gamma=2.5,3,4,5,7,10,\infty$. Notice that, for ease of comparison, the values
associated with $\gamma=\infty$ in the rightmost column of Table 2 have been multiplied by 100, since these values are quite close to 0.
The tables presented were obtained using the R package $stargazer$, see \cite{Hla2015}.
From Table 2, we see that the speed of convergence to the asymptotic null distribution is slow and depends on the value of $\gamma$ as well as on the dimension of the data.

% For tables use
\begin{table}
%\label{tab1}
%\caption{Critical values for ($16 \gamma^{2+d/2}/ \pi^{d/2}) T_{n,\gamma}$, $\alpha = 0.05$}
%\label{tab1}       % Give a unique label
% For LaTeX tables use
\renewcommand{\arraystretch}{0.75}
\begin{tabular}{|p{0.7cm}|ccccccc|}
\hline
\hline
    &       &     &     & $\gamma$ &       &      &   \\
$n$ & $2.5$ & $3$ & $4$ & $5$      & $7$   & $10$ & $\infty$ \\
\hline
&  & & & $d=1$ & & & \\
\hline
20  & $120.42$ & $105.16$ & $88.41$ & $79.48$ & $70.66$ & $64.74$ & $265.14$ \\
 50  & $219.25$ & $173.63$ & $130.90$ & $111.24$ & $93.12$ & $82.02$ & $125.20$ \\
 100  & $294.01$ & $218.62$ & $154.20$ & $126.76$ & $102.94$ & $89.04$ & $66.13$ \\
 200  & $361.62$ & $254.88$ & $170.48$ & $136.65$ & $108.60$ & $92.88$ & $33.89$ \\
 300  & $395.67$ & $271.05$ & $176.56$ & $140.20$ & $110.50$ & $94.19$ & $22.79$ \\
 \hline  & & & & $d=2$ & & & \\
\hline
20  & $535.50$ & $413.65$ & $300.96$ & $249.60$ & $203.09$ & $174.73$ & $628.97$ \\
 50  & $1,086.28$ & $737.69$ & $464.16$ & $356.56$ & $268.57$ & $220.27$ & $291.96$ \\
 100 & $1,516.50$ & $947.61$ & $546.91$ & $402.84$ & $292.63$ & $235.24$ & $152.02$ \\
 200  & $1,867.44$ & $1,089.76$ & $585.61$ & $419.94$ & $299.04$ & $238.36$ & $77.02$ \\
 300  & $2,035.48$ & $1,141.98$ & $595.36$ & $422.25$ & $299.14$ & $238.42$ & $51.52$ \\
 \hline &  & & & $d=3$ & & & \\
\hline
20  & $1,460.39$ & $1,044.12$ & $695.34$ & $549.27$ & $423.00$ & $350.28$ & $1,157.20$ \\
 50  & $3,444.99$ & $2,095.29$ & $1,162.22$ & $831.36$ & $580.61$ & $451.38$ & $537.68$ \\
 100  & $5,054.15$ & $2,781.23$ & $1,384.12$ & $941.78$ & $628.33$ & $477.34$ & $278.23$ \\
 200  & $6,463.29$ & $3,267.95$ & $1,495.75$ & $980.63$ & $638.18$ & $481.29$ & $140.44$ \\
 300  & $7,108.14$ & $3,439.05$ & $1,508.96$ & $977.42$ & $633.30$ & $478.42$ & $93.78$ \\
 \hline &  & & & $d=5$ & & & \\
\hline
20  & $6,346.44$ & $4,065.35$ & $2,389.07$ & $1,759.71$ & $1,257.17$ & $986.77$ & $2,903.55$ \\
 50  & $20,164.36$ & $10,268.49$ & $4,655.52$ & $2,988.80$ & $1,862.85$ & $1,340.01$ & $1,361.04$ \\
 100  & $34,187.51$ & $15,193.42$ & $5,934.77$ & $3,545.86$ & $2,070.71$ & $1,439.25$ & $705.30$ \\
 200  & $47,436.90$ & $18,844.25$ & $6,578.37$ & $3,746.81$ & $2,114.14$ & $1,450.44$ & $355.91$ \\
 300  & $54,128.44$ & $20,321.98$ & $6,715.41$ & $3,749.27$ & $2,091.10$ & $1,436.11$ & $237.36$ \\
 \hline
\end{tabular}
\begin{center}
{\bf Table 2:}  Critical values for ($16 \gamma^{2+d/2}/ \pi^{d/2}) T_{n,\gamma}$, $\alpha = 0.05$
\end{center}
\end{table}

\subsection{Power results}

In the power study presented in this section we simulate from several univariate and multivariate distributions. The critical values presented in Table 2
 are used in order to calculate the power results for the newly proposed test. For the remaining tests, the critical values used are also calculated based on $10^6$ simulations under $H_0$. The power estimates presented below are based on 10 000 random samples. Notice that, due to the computational complexity of $HM_n$, the number of simulations used in order to obtain both critical values and power estimates for this test are reduced by a factor of 10. This step is necessary in order to ensure that the numerical results can be obtained within a reasonable time.

Power results are reported for a sample size of $n=50$ and data dimensions $d=1,2,3,5$. A nominal significance level of $\alpha = 0.05$ is used throughout. As is pointed out in \cite{MM2005}, the maximum possible standard error for each power estimate is $0.005$.
Thus, the 99\% confidence interval for each reported power estimate is contained in the interval obtained by subtracting 1\% from, and adding 1\% to, the stated power.

The results for the univariate case are given in Table 3. The entries are percentages of rejection of $H_0$, rounded to the nearest integer.
The power against the standard normal distribution shows that the nominal level is maintained very closely.
As for the alternative distributions, NMIX1 and NMIX2 denote mixtures of the normal distributions $\textrm{N}(0,1)$ and $\textrm{N}(0,4)$.
The mixture NMIX1 gives equal weight to these distributions, while NMIX2 is obtained
when the probability of sampling from the standard normal distribution is 0.75. The remaining alternative distributions considered are $t$-distributions with 3, 5 and 10 degress of freedom, the lognormal distributions with parameters $(0,1/2)$ and $(0,1/4)$ (denoted by $LN(\cdot)$), the $\chi^2$-distributions with 5 and 15 degrees of freedom, the standard logistic distribution, the Weibull distributions with shape parameters 10 and 20, the Pearson type VII distributions with 5 and 10 degrees of freedom (denoted by $P_{VII}(\cdot)$), and the skew-normal law with skewness parameters 3 and 5 (denoted by $SN(\cdot)$), see \cite{Azz1985}.

The R Package $PearsonDS$ contains the function $rpearsonVII$, which can be employed to simulate random variables from this distribution,
see \cite{BK2017}. The R Package $sn$ disposes of the function $rsn$, which can be used to simulate random variates from the skew normal distribution, see \cite{Azz2017}.

Tables 3,4, 5 and 6  report the powers calculated in the cases where $d$ equals 1, 2, 3 and 5 respectively.
Note that the subscript $n$ in the names of the test statistics is omitted in the tables in order to save space.

\begin{table}[!htbp] \centering
  %\label{tab2}
  \renewcommand{\arraystretch}{0.75}
\begin{tabular}{@{\extracolsep{5pt}} p{2.5cm}p{0.6cm}p{0.6cm}p{0.6cm}p{0.6cm}p{0.6cm}p{0.6cm}p{0.6cm}p{0.6cm}p{0.6cm}p{0.6cm}}
\hline
\hline
&  $CvM$ & $AD$ & $SW$ & $JB$ & $Z_{3}$ & $Z_{15}$ & $T_{2.5}$ & $T_{5}$ & $T_{10}$ & $T_{\infty}$ \\
\hline
$\textrm{N}(0,1)$   & $5$ & $5$ & $5$ & $5$ & $5$ & $5$ & $5$ & $5$ & $5$ & $5$ \\
\textrm{NMIX1}       & $18$ & $20$ & $21$ & $28$ & $24$ & $20$ & $24$ & $23$ & $21$ & $18$ \\
\textrm{NMIX2}       & $19$ & $22$ & $28$ & $37$ & $34$ & $28$ & $34$ & $32$ & $30$ & $26$ \\
$t(3)$                       & $58$ & $61$ & $63$ & $69$ & $65$ & $57$ & $65$ & $63$ & $60$ & $52$ \\
$t(5)$                       & $28$ & $31$ & $37$ & $44$ & $41$ & $35$ & $41$ & $39$ & $37$ & $32$ \\
$t(10)$                     & $12$ & $13$ & $16$ & $21$ & $20$ & $17$ & $20$ & $20$ & $19$ & $17$ \\
$LN(0,\frac{1}{2})$ & $83$ & $87$ & $93$ & $85$ & $80$ & $89$ & $76$ & $85$ & $88$ & $91$ \\
$LN(0,\frac{1}{4})$ & $31$ & $35$ & $44$ & $39$ & $37$ & $45$ & $34$ & $40$ & $44$ & $47$ \\
$\chi^2(5)$               & $74$ & $81$ & $89$ & $75$ & $69$ & $82$ & $62$ & $74$ & $80$ & $83$ \\
$\chi^2(15)$             & $30$ & $34$ & $43$ & $36$ & $34$ & $43$ & $31$ & $38$ & $41$ & $45$ \\
$Logistic(0,1)$          & $14$ & $16$ & $19$ & $26$ & $24$ & $20$ & $24$ & $23$ & $21$ & $19$ \\
$Weibull(10)$            & $25$ & $28$ & $34$ & $30$ & $28$ & $36$ & $25$ & $31$ & $35$ & $37$ \\
$Weibull(20)$            & $39$ & $44$ & $53$ & $46$ & $44$ & $53$ & $40$ & $48$ & $52$ & $55$ \\
$P_{VII}(5)$             & $27$ & $30$ & $36$ & $43$ & $40$ & $35$ & $41$ & $39$ & $37$ & $32$ \\
$P_{VII}(10)$           & $10$ & $12$ & $16$ & $21$ & $20$ & $17$ & $20$ & $19$ & $18$ & $16$ \\
$SN(3)$                     & $31$ & $34$ & $40$ & $32$ & $30$ & $39$ & $25$ & $33$ & $37$ & $41$ \\
$SN(5)$                     & $53$ & $59$ & $67$ & $49$ & $43$ & $58$ & $36$ & $49$ & $55$ & $61$ \\
\hline
\end{tabular}
\begin{center}
{\bf Table 3:}  Monte Carlo power estimates in the univariate case, $\alpha = 0.05$
\end{center}
    \end{table}

The results shown in Table 3 indicate that the newly proposed class of tests exhibit substantial power against the distributions considered.
$T_\infty$ outperforms each of the competing tests in terms of the estimated powers against the $LN(0,1/4)$, $\chi^2(15)$ and $SN(3)$ laws,
as well as both of the Weibull distributions considered. The newly proposed tests also proves to be serious competitors against
each of the remaining distributions considered, especially for small values of $\gamma$ (in which case $T_{n,\gamma}$ is often only outperformed by the Jarque-Bera test).

We now turn our attention to the case $d>1$. As was the case for the univariate tests, the powers against 16 alternative
distributions are reported for each of the data dimensions considered. The powers of each of the tests against the standard normal
distribution are also included in the relevant tables. The alternative distributions considered include mixtures of normal laws,
distributions with independent marginals, distributions where each marginal is normal except for one, a spherically symmetric distribution,
 and a distribution for which every marginal is normal, but the joint distribution does not follow the normal law.

The parameter combinations used for the mixtures of normal distributions were taken from \cite{SR2005}.
Let $p\textrm{N}_d(\mu_1,\Sigma_1)+(1-p)\textrm{N}_d(\mu_2,\Sigma_2)$
denote a normal mixture, where the probability of sampling from $\textrm{N}_d(\mu_1,\Sigma_1)$ is $p$ and the probability of sampling from $\textrm{N}_d(\mu_2,\Sigma_2)$ is $1-p$.
Let $\mu=0$ and $\mu=3$ denote $d$-dimensional column vectors of 0's and 3's, respectively, and let $B_d$ denote a ($d \times d$)-matrix containing 1's on the main diagonal
and 0.9's on each off diagonal entry. The normal mixtures are constructed by combining $\textrm{N}_d(0,\textrm{I}_d)$, $\textrm{N}_d(3,\textrm{I}_d)$ and $\textrm{N}_d(0,B_d)$.
The first mixture, denoted by \textrm{NMIX1}, is $0.9\textrm{N}_d(0,\textrm{I}_d)+0.1\textrm{N}_d(3,\textrm{I}_d)$. This distribution is skewed with heavy tails. The second mixture,
denoted by \textrm{NMIX2}, is $0.9\textrm{N}_d(0,B_d)+0.1\textrm{N}_d(0,\textrm{I}_d)$. This is a symmetric, heavy-tailed distribution.
In addition, we included two multivariate $t$-distributions; the $t_\nu(0,\textrm{I}_d)$-distribution for $\nu=5$ and $\nu=10$. Next, we included distributions
with independent marginals, the latter being the $\chi^2$-distribution with 15 and 20 degrees of freedom respectively, the $logistic(0,1)$ distribution, the gamma distributions with parameters $(5,1)$ and $(4,2)$,
as well as the Pearson Type VII distributions with 10 and 20 degrees of freedom.

Three $d$-dimensional distributions are obtained by combining $d-1$ independent standard normal marginals with one non-normal distribution.
This distribution is denoted by $\textrm{N}(0,1)^{d-1} \otimes F$, where $F$ denotes the non-normal marginal distribution. The three alternatives considered for $F$ are
the $\chi^2$-distributions with 5  and 10 degrees of freedom, respectively, as well as the $t$-distribution with 3 degrees of freedom.

Spherically symmetric distributions can be defined in  R using the $EllipticalDistribution$ function from the R package $distrEllipse$, see \cite{RKSC2006}.
Tables 4--6 display the estimated powers of the various tests considered against the $d$-dimensional spherically symmetric distribution,
where the radius of the distribution follows a lognormal distribution with parameters 0 and 0.5. This distribution is denoted by $\mathcal{S}^d(LN(0,1/2))$.

Let $\rho_d$ and $\rho_d'$ denote positive definite ($d\times d$)-matrices with 1's on the main diagonal, where $\rho_d$ has the constant $\rho$ and $\rho_d'$ the constant $-\rho$ on each off diagonal entry.
 The final distribution considered is the mixture $0.5 \textrm{N}_d(0,\rho_d) + 0.5 \textrm{N}_d(0,\rho_d')$. This distribution is a non-normal $d$-variate distribution with normal marginals.

\begin{table}[!htbp] \centering
  %\label{tab3}
  \renewcommand{\arraystretch}{0.75}
\begin{tabular}{@{\extracolsep{5pt}} p{2.5cm}p{0.6cm}p{0.6cm}p{0.6cm}p{0.6cm}p{0.6cm}p{0.6cm}p{0.6cm}p{0.6cm}p{0.6cm}p{0.6cm}}
\hline
\hline
& $MS$ & $MK$ & $HZ$ & $EN$ & $HM$ & $HJ$ & $T_{2.5}$ & $T_{5}$ & $T_{10}$ & $T_{\infty}$ \\
\hline
$\textrm{N}(0,1)$                                  & $5$ & $5$ & $5$ & $5$ & $5$ & $5$ & $5$ & $5$ & $5$ & $5$ \\
$\textrm{NMIX1}$                                  & $85$ & $34$ & $75$ & $82$ & $57$ & $73$ & $48$ & $69$ & $80$ & $86$ \\
$\textrm{NMIX2}$                                  & $44$ & $48$ & $29$ & $38$ & $57$ & $53$ & $55$ & $54$ & $52$ & $44$ \\
$\textrm{t}_{5}(0,\textrm{I}_2)$         & $53$ & $62$ & $42$ & $51$ & $67$ & $60$ & $60$ & $60$ & $58$ & $53$ \\
$\textrm{t}_{10}(0,\textrm{I}_2)$       & $24$ & $26$ & $14$ & $19$ & $32$ & $29$ & $29$ & $29$ & $28$ & $25$ \\
$(\chi^2(15))^2$                                    & $49$ & $19$ & $34$ & $42$ & $26$ & $41$ & $30$ & $39$ & $45$ & $52$ \\
$(\chi^2(20))^2$                                    & $40$ & $16$ & $27$ & $33$ & $24$ & $34$ & $25$ & $32$ & $37$ & $42$ \\
$\textrm{Logistic}(0,1)^2$                     & $24$ & $27$ & $15$ & $19$ & $33$ & $28$ & $28$ & $29$ & $28$ & $25$ \\
$\textrm{Gamma}(5,1)^2$                     & $67$ & $27$ & $52$ & $61$ & $38$ & $57$ & $41$ & $54$ & $62$ & $70$ \\
$\textrm{Gamma}(4,2)^2$                     & $76$ & $32$ & $64$ & $72$ & $42$ & $66$ & $48$ & $62$ & $71$ & $78$ \\
$P_{VII}(10)^2$                                     & $20$ & $21$ & $11$ & $14$ & $27$ & $23$ & $24$ & $24$ & $23$ & $20$ \\
$P_{VII}(20)^2$                                     & $11$ & $10$ & $7$   & $8$   & $14$ & $12$ & $13$ & $13$ & $12$ & $12$ \\
$\textrm{N}(0,1)\otimes$$t(3)$              & $47$ & $52$ & $42$ & $49$ & $61$ & $55$ & $56$ & $56$ & $54$ & $47$ \\
$\textrm{N}(0,1)\otimes$$\chi^2(5)$     & $63$ & $25$ & $52$ & $60$ & $36$ & $52$ & $39$ & $49$ & $57$ & $65$ \\
$\textrm{N}(0,1)\otimes$$\chi^2(10)$   & $38$ & $15$ & $26$ & $32$ & $21$ & $32$ & $24$ & $30$ & $35$ & $40$ \\
$\mathcal{S}^2(LN(0,\frac{1}{2}))$       & $26$ & $25$ & $15$ & $21$ & $29$ & $30$ & $31$ & $31$ & $29$ & $26$ \\
$\textrm{NM}_2(\rho =0.2)$                   & $6$ & $6$ & $5$ & $6$ & $6$ & $6$ & $6$ & $6$ & $6$ & $6$ \\
\hline
\end{tabular}
\begin{center}
{\bf Table 4:} Monte Carlo power estimates in the multivariate case for $d=2$
\end{center}
\end{table}

\begin{table}[!htbp] \centering
%  \caption{Monte Carlo power estimates in the multivariate case for $d=3$}
  %\label{tab4}
  \renewcommand{\arraystretch}{0.75}
\begin{tabular}{@{\extracolsep{5pt}} p{2.5cm}p{0.6cm}p{0.6cm}p{0.6cm}p{0.6cm}p{0.6cm}p{0.6cm}p{0.6cm}p{0.6cm}p{0.6cm}p{0.6cm}}
\hline
\hline
& $MS$ & $MK$ & $HZ$ & $EN$ & $HM$ & $HJ$ & $T_{2.5}$ & $T_{5}$ & $T_{10}$ & $T_{\infty}$ \\
\hline
$\textrm{N}(0,1)^3$                                   & $5$ & $5$ & $5$ & $5$ & $5$ & $5$ & $5$ & $5$ & $5$ & $5$ \\
$\textrm{NMIX1}$                                        & $89$ & $36$ & $81$ & $91$ & $59$ & $72$ & $43$ & $66$ & $82$ & $91$ \\
$\textrm{NMIX2}$                                        & $71$ & $76$ & $49$ & $66$ & $79$ & $79$ & $79$ & $80$ & $78$ & $72$ \\
$\textrm{t}_{5}(0,\textrm{I}_3)$               & $68$ & $78$ & $55$ & $68$ & $77$ & $73$ & $71$ & $73$ & $73$ & $69$ \\
$\textrm{t}_{10}(0,\textrm{I}_3)$             & $34$ & $38$ & $18$ & $27$ & $35$ & $38$ & $36$ & $38$ & $38$ & $34$ \\
$(\chi^2(15))^3$                                         & $52$ & $21$ & $35$ & $49$ & $27$ & $42$ & $31$ & $39$ & $47$ & $55$ \\
$(\chi^2(20))^3$                                         & $40$ & $16$ & $26$ & $37$ & $21$ & $33$ & $24$ & $30$ & $36$ & $44$ \\
$\textrm{Logistic}(0,1)^3$                          & $28$ & $30$ & $15$ & $22$ & $33$ & $31$ & $30$ & $31$ & $31$ & $28$ \\
$\textrm{Gamma}(5,1)^3$                          & $72$ & $30$ & $53$ & $69$ & $39$ & $58$ & $41$ & $53$ & $65$ & $75$ \\
$\textrm{Gamma}(4,2)^3$                          & $80$ & $36$ & $65$ & $79$ & $46$ & $66$ & $47$ & $61$ & $73$ & $83$ \\
$P_{VII}(10)^3$                                          & $22$ & $22$ & $10$ & $16$ & $24$ & $25$ & $25$ & $26$ & $25$ & $23$ \\
$P_{VII}(20)^3$                                          & $12$ & $10$ & $6$   & $8$   & $14$ & $13$ & $13$ & $13$ & $13$ & $12$ \\
$\textrm{N}(0,1)^2\otimes$$t(3)$              & $42$ & $43$ & $29$ & $40$ & $54$ & $48$ & $49$ & $49$ & $48$ & $43$ \\
$\textrm{N}(0,1)^2\otimes$$\chi^2(5)$     & $47$ & $18$ & $33$ & $46$ & $28$ & $39$ & $29$ & $36$ & $43$ & $51$ \\
$\textrm{N}(0,1)^2\otimes$$\chi^2(10)$   & $26$ & $12$ & $17$ & $24$ & $16$ & $22$ & $17$ & $21$ & $24$ & $28$ \\
$\mathcal{S}^3(LN(0,\frac{1}{2}))$           & $53$ & $58$ & $18$ & $43$ & $62$ & $58$ & $57$ & $58$ & $58$ & $54$ \\
$\textrm{NM}_3(\rho =0.2)$                       & $8$ & $7$ & $5$ & $6$ & $7$ & $8$ & $8$ & $8$ & $8$ & $8$ \\
\hline
\end{tabular}
\begin{center}
{\bf Table 5:} Monte Carlo power estimates in the multivariate case for $d=3$
\end{center}
\end{table}

\begin{table}[!htbp] \centering
  %\label{tab5}
  \renewcommand{\arraystretch}{0.75}
\begin{tabular}{@{\extracolsep{5pt}} p{2.5cm}p{0.6cm}p{0.6cm}p{0.6cm}p{0.6cm}p{0.6cm}p{0.6cm}p{0.6cm}p{0.6cm}p{0.6cm}p{0.6cm}}
\hline
\hline
& $MS$ & $MK$ & $HZ$ & $EN$ & $HM$ & $HJ$ & $T_{2.5}$ & $T_{5}$ & $T_{10}$ & $T_{\infty}$ \\
\hline
$\textrm{N}(0,1)^5$                                  & $5$ & $5$ & $5$ & $5$ & $5$ & $5$ & $5$ & $5$ & $5$ & $5$ \\
$\textrm{NMIX1}$                                       & $82$ & $33$ & $74$ & $94$ & $43$ & $58$ & $34$ & $51$ & $68$ & $86$ \\
$\textrm{NMIX2}$                                       & $94$ & $94$ & $68$ & $89$ & $95$ & $95$ & $95$ & $96$ & $95$ & $94$ \\
$\textrm{t}_{5}(0,\textrm{I}_5)$              & $88$ & $94$ & $72$ & $88$ & $89$ & $90$ & $86$ & $89$ & $90$ & $89$ \\
$\textrm{t}_{10}(0,\textrm{I}_5)$            & $54$ & $58$ & $23$ & $45$ & $51$ & $55$ & $51$ & $55$ & $57$ & $55$ \\
$\chi^2(15)^5$                                         & $51$ & $22$ & $30$ & $52$ & $26$ & $39$ & $29$ & $36$ & $44$ & $56$ \\
$\chi^2(20)^5$                                         & $39$ & $16$ & $22$ & $39$ & $20$ & $30$ & $23$ & $28$ & $33$ & $42$ \\
$\textrm{Logistic}(0,1)^5$                          & $33$ & $34$ & $13$ & $25$ & $31$ & $34$ & $31$ & $34$ & $36$ & $33$ \\
$\textrm{Gamma}(5,1)^5$                          & $72$ & $33$ & $49$ & $74$ & $37$ & $55$ & $40$ & $51$ & $63$ & $76$ \\
$\textrm{Gamma}(4,2)^5$                          & $81$ & $40$ & $60$ & $84$ & $40$ & $64$ & $47$ & $59$ & $72$ & $85$ \\
$P_{VII}(10)^5$                                        & $27$ & $25$ & $9$   & $19$ & $26$ & $28$ & $26$ & $28$ & $29$ & $27$ \\
$P_{VII}(20)^5$                                        & $12$ & $9$   & $6$   & $9$   & $12$ & $12$ & $11$ & $12$ & $12$ & $12$ \\
$\textrm{N}(0,1)^4\otimes$$t(3)$              & $35$ & $32$ & $16$ & $30$ & $42$ & $39$ & $38$ & $39$ & $39$ & $35$ \\
$\textrm{N}(0,1)^4\otimes$$\chi^2(5)$     & $28$ & $13$ & $16$ & $28$ & $19$ & $23$ & $19$ & $22$ & $25$ & $31$ \\
$\textrm{N}(0,1)^4\otimes$$\chi^2(10)$   & $16$ & $8$   & $10$ & $15$ & $13$ & $14$ & $12$ & $13$ & $15$ & $18$ \\
$\mathcal{S}^5(LN(0,\frac{1}{2}))$           & $89$ & $95$ & $77$ & $90$ & $90$ & $90$ & $86$ & $90$ & $91$ & $89$ \\
$\textrm{NM}_5(\rho =0.2)$                       & $12$ & $9$   & $5$   & $11$ & $12$ & $13$ & $12$ & $13$ & $13$ & $12$ \\
\hline
\end{tabular}
\begin{center}
{\bf Table 6:} Monte Carlo power estimates in the multivariate case for $d=5$
\end{center}
\end{table}

As was the case in the univariate setting, the newly proposed test is associated with several high powers reported in Tables 4, 5, and 6.
When comparing the results for the distributions with independent marginals, we see that $T_{\infty}$ outperforms each of the competitors against both of
the distributions with $\chi^2$ marginals as well as both of the distributions with gamma marginals considered. This is also the case aginst
the $\textrm{N}(0,1)^{d-1}\otimes$$\chi^2(5)$ and $\textrm{N}(0,1)^{d-1}\otimes$$\chi^2(10)$ distributions. The mentioned predominance is for each of the
 data dimensions considered. Furthermore, the newly proposed test statistic shows high power against the remaining distributions for finite values of $\gamma$.
  Specifically, when $d=2$ the new test outperforms its competitors against the spherically symmetric distribution. In the case $d=3$, none of the
   competing tests are able to outperform the newly proposed class of tests against the second normal mixture considered, the $t$-distribution with 10 degrees
    of freedom or the $PVII(10)$ distribution. The corresponding list of distributions in the case where $d=5$ is obtained by substituting the $t$-distribution
     with 10 degrees of freedom for the $t$-distribution with 5 degrees of freedom and adding the $PVII(20)$ distribution. Finally, none of the competing tests
     is able to outperform the newly proposed class of tests against the $NM_d(\rho=0.2)$ distribution.

In most of the cases considered, the power of the newly proposed class of tests is a monotone function of $\gamma$. Based on the numerical results presented, it is recommeded that $\gamma=5$ be used when performing the test as this value results in reasonably high power against the majority of the alternatve distributions considered.

%%%%%%%%%%%%%%%%%%%%%%%%%%%%%%%%%%%%%%%%%%%%%%%%%%%%%%%%%%%%%%%%%%%%%%%%%%%%%%%%%%%%%%%%%%%%%%%%%%%%%%%%%%%%%%%%%%%%%%%%%%%%%%%%%%%%%%%%%%%%%%%%%%%%%%%%%%%%%%%%%%%%%%%%%%%%%%
%
%
%                   A real data example
%
%
%%%%%%%%%%%%%%%%%%%%%%%%%%%%%%%%%%%%%%%%%%%%%%%%%%%%%%%%%%%%%%%%%%%%%%%%%%%%%%%%%%%%%%%%%%%%%%%%%%%%%%%%%%%%%%%%%%%%%%%%%%%%%%%%%%%%%%%%%%%%%%%%%%%%%%%%%%%%%%%%%%%%%%%%%%%%%%

\section{A real data example}\label{secrealdata}

The payoff function of certain types of financial derivatives depends on the joint behaviour of multiple stocks or indexes; an important example is the class of basket options. When calculating the price of a basket option, it is often assumed that the log-returns of the stocks or indexes considered are realized from a multivariate normal distribution (this assumption is an extention of the celebrated Black-Merton-Scholes model for options on a single stock or index). As a result, testing the hypothesis that observed financial log-returns follow a multivariate normal law is of interest when pricing basket options. For more details regarding the pricing of these options, the interested reader is referred to \cite{CF2016}.

As a practical application, we consider the log-returns associated with three major indexes traded in the financial market of the United States. 50 daily log-returns were calculated for the period ending 29 December 2017, the relevant prices were downloaded from http://finance.yahoo.com. The first index considered is the Dow Jones Industrial Average (DJIA) index; this index is comprised of a price-weighted average of 30 large publicly owned companies. The second is the Standard \& Poor 500 (S\&P 500); a market-capitalization weighted index comprising 500 large companies. Finally, we consider the log-returns of the National Association of Securities Dealers Automated Quotations (NASDAQ) composite index. We are interested in testing the hypothesis that the log-returns are realized from a  multivariate normal distribution using the newly proposed test.

Table 7} shows the estimated p-values associated with the newly proposed tests for various values of $\gamma$, the reported p-values were obtained using one million Monte Carlo simulations in each case. The results indicate that the hypothesis of multivariate normality is rejected at a 1\% significance level for each value of $\gamma$ considered.

\begin{table}[!htbp] \centering
    %\label{tab6}
    \renewcommand{\arraystretch}{0.75}
\begin{tabular}{ccccccc}
\hline
\hline
$\gamma=2.5$ & $\gamma=3$ & $\gamma=4$ & $\gamma=5$ & $\gamma=7$ & $\gamma=10$ & $\gamma=\infty$ \\
\hline \\[-1.8ex]
$0.0072$ & $0.0061$ & $0.0044$ & $0.0033$ & $0.00210$ & $0.0013$ & $0.0002$\\
\hline \\[-1.8ex]
\end{tabular}
\begin{center}
{\bf Table 7:} p-values associated with the newly proposed tests
\end{center}
\end{table}

\section{Conclusion}\label{secconclus}

We proposed and studied a new class of affine invariant tests for normality in any dimension that are based on a partial differential
equation involving the moment generating function. Some properties of the limit null distribution of the test statistic  $T_{n,\gamma}$ have been derived,
and the consistency of this class of tests against general alternatives has been proved under some mild conditions.
For fixed $n$, the test statistic  $T_{n,\gamma}$, after suitable scaling, approaches a linear combination of two  measures of multivariate skewness as $\gamma \to \infty$.

A Monte Carlo study investigates the finite-sample performance of $T_{n,\gamma}$  compared to those of competing tests in the univariate
and multivariate settings. The competing tests considered for univariate normality comprise four well-known tests, while,
in the multivariate case, we include four prominent classical tests for multinormality and two very recent tests.
The numerical results indicate that $T_{n,\gamma}$ often exhibits power greater than those associated with several
of its competitors, both in univariate and multivariate settings. Based on the numerical results obtained, it is recommended that $\gamma=5$ is used when performing the test.\\[1mm]

{\bf Acknowledgement:} The authors thank Tobias Jahnke for providing the proof of Proposition  \ref{proplinalgebra}.
The second author's work is based on research supported by the National Research Foundation, South Africa (Research chair: Non-parametric, Robust Statistical Inference and Statistical Process Control, Grant number 71199). Opinions expressed and conclusions arrived at are those of the authors and are not necessarily to be attributed to the NRF.

\section{Appendix: Technical proofs} \label{sectechnical}

\begin{prop}\label{propcnt} We have $\|C_n\|_{\oplus} = o_\mathbb{P}(1)$, where $C_n$ is given in (\ref{cnt}).
\end{prop}

\noindent {\sc Proof.}
Putting
$
\xi_j(t) := \exp(t^\top X_j \left(t^\top \Delta_{n,j}\right)^2 \exp\left(\Theta_{n,j} t^\top \Delta_{n,j}\right)/2
$,
we have
\begin{eqnarray*}
\|C_n(t)\|^2 & = & \frac{1}{n} \sum_{j,k=1}^n \xi_j(t) \xi_k(t) \left(X_j-t+ \Delta_{n,j}\right)^\top \left(X_k-t+ \Delta_{n,k}\right).
\end{eqnarray*}
Recall $F_n$ from (\ref{defmaxnorm}) and put $\Lambda_n := \|S_n^{-1/2}-\textrm{I}_d\|_2 F_n + \|S_n^{-1/2}\|_2 \|\overline{X}_n\|$. We have
$|\Theta_{n,j} t^\top \Delta_{n,j}| \le \|t\|\, \Lambda_n$. Furthermore, using $|t^\top (X_j+X_k)| \le 2 \, \|t\| \, F_n$ and
\[
\left(t^\top \Delta_{n,j}\right)^2 \le 2 \|t\|^2 \, \|S_n^{-1/2}- \textrm{I}_d\|^2_2 \, \|X_j\|^2 + 2 \|t\|^2 \, \|S_n^{-1/2}\|^2_2 \, \|\overline{X}_n\|^2,
\]
as well as
\[
%\begin{eqnarray*}
|\left(X_j\! -\! t\! +\! \Delta_{n,j}\right)^\top \! \left(X_k\! -\! t \! + \!  \Delta_{n,k}\right)| \le %\qquad \qquad \qquad \qquad \\
%\qquad \qquad \qquad \qquad
\left(\|X_j\|\! + \!  \|t\| \! + \!  \|\Delta_{n,j}\|\right) \left(\|X_k\|\!  + \!  \|t\| \!  + \!  \|\Delta_{n,k}\|\right),
%\end{eqnarray*}
\]
we obtain
\begin{eqnarray*}
\|C_n(t)\|^2 & \le & e^{2 \|t\| \Gamma_n} \frac{1}{n}\sum_{j,k=1}^n \left( \Big{\{} \|t\|^2 \, \|S_n^{-1/2}-\textrm{I}_d\|^2_2 \, \|X_j\|^2 + \|t\|^2 \, \|S_n^{-1/2}\|^2_2 \, \|\overline{X}_n\|^2\Big{\}} \right.  \\
& & \qquad \qquad \qquad \times \Big{\{} \|t\|^2 \, \|S_n^{-1/2}-\textrm{I}_d\|^2_2 \, \|X_k\|^2 + \|t\|^2 \, \|S_n^{-1/2}\|^2_2 \, \|\overline{X}_n\|^2\Big{\}} \\
& & \qquad \qquad \qquad \left. \times \Big{\{} \|X_j\| + \|t\| + \|\Delta_{n,j}\| \Big{\}} \Big{\{} \|X_k\| + \|t\| + \|\Delta_{n,k}\|\Big{\}}\right),
\end{eqnarray*}
where $\Gamma_n = F_n + \Lambda_n$. Expanding the curly brackets, the leading terms are those that do not involve any of $\Delta_{n,j}$. We concentrate on
\[
S_n(t) := \textrm{e}^{2 \|t\| \Gamma_n} \frac{1}{n}\sum_{j,k=1}^n \left(\|t\|^4  \,  \|S_n^{-1/2}-\textrm{I}_d\|^4_2 \, \|X_j\|^3 \, \|X_k\|^3\right),
\]
which originates from choosing the first term within each of the curly brackets. The other terms are treated similarly.  Notice that
\[
\int_{\mathbb{R}^d} S_n(t) \, w_\gamma(t) \, \textrm{d}t = \|S_n^{-1/2}-\textrm{I}_d\|^4_2 \cdot n \cdot  \Big{(}\frac{1}{n}\sum_{j=1}^n \|X_j\|^3\Big{)}^2 \cdot
\int_{\mathbb{R}^d} \|t\|^4 \textrm{e}^{2\|t\|\Gamma_n} \, w_\gamma(t) \, \textrm{d} t,
\]
and that $\|S_n^{-1/2}-\textrm{I}_d\|^2_2 = O_\mathbb{P}(n^{-2})$, $n^{-1}\sum_{j=1}^n \|X_j\| = O_\mathbb{P}(1)$.
By Proposition 10.2 of \cite{HJM2017}, the integral is of order $O_\mathbb{P}(\Gamma_n^{d+3}) \exp(\Gamma_n^2/\gamma)$ (notice that
$1+\gamma$ in that paper corresponds to (our) $\gamma$). From display (10.6) and display (10.7) of \cite{HJM2017} we have
$
\Gamma_n = O_\mathbb{P}(\sqrt{\log n})$ and $\exp\left( \Gamma_n^2/\gamma\right) = n^{2/\gamma} \cdot \left(\log n\right)^{(d-2)/\gamma} \, O_\mathbb{P}(1)$.
Since $\gamma > 2$, it follows that
\[
\int_{\mathbb{R}^d} S_n(t) \, w_\gamma(t) \, \textrm{d}t = O_\mathbb{P}\left(n^{\frac{2}{\gamma}-1}\right) O_\mathbb{P}\left(\left(\log n\right)^{\frac{d+3}{2} + \frac{d-2}{\gamma}} \right) \ = \ o_\mathbb{P}(1). \quad \bewend
\]

\begin{prop}\label{propbnt} For $B_n$ given in (\ref{bnt}), we have
\[
B_n(t) = - \frac{1}{\sqrt{n}} \, \exp \left(\frac{\|t\|^2}{2}\right) \sum_{j=1}^n \left( \frac{X_jX_j^\top t - t}{2} + X_j \right) + o_\mathbb{P}(1).
\]
\end{prop}

\noindent{\sc Proof.} Observe that $B_n(t) = n^{-1/2} \sum_{j=1}^n \exp(t^\top X_j) \Delta_{n,j} + R_n(t)$, where
\[
R_n(t) = \frac{1}{\sqrt{n}} \sum_{j=1}^n \textrm{e}^{t^\top X_j} (t^\top \Delta_{n,j} ) \, \Delta_{n,j}.
\]
Use
\[
\|R_n(t)\|^2 \le \Big{(} \frac{1}{n}\sum_{j=1}^n \textrm{e}^{t^\top X_j} \Big{)}^2 \, n \, \|t\|^2  \, \max_{i=1,\ldots,n} \|\Delta_{n,i}\|^4
\]
and
$
\|\Delta_{n,j} \| \le \|S_n^{-1/2}-\textrm{I}_d\|_2 \, \|F_n\| + \|S_n^{-1/2}\| \, \|\overline{X}_n\|
$
with $F_n$ given in (\ref{defmaxnorm}) together with $F_n= O_\textrm{P}(\sqrt{\log n})$ (see Prop. 10.1. of \cite{HJM2017}) and
\begin{equation}\label{eabsch}
\mathbb{E}\Big{[}\Big{(}\frac{1}{n}\sum_{j=1}^n \textrm{e}^{t^\top X_j} \Big{)}^2\Big{]} = \mathbb{E}\Big{[}\frac{1}{n^2}\sum_{j,k=1}^n \textrm{e}^{t^\top(X_j+X_k)}\Big{]}
\le \textrm{e}^{2\|t\|^2} + \textrm{e}^{\|t\|^2}
\end{equation}
to show  $R_n = o_\mathbb{P}(1)$. Next,
$
B_n(t) - R_n(t) = B_{n,1}(t) - B_{n,2}(t) - B_{n,3}(t),
$
where
\begin{eqnarray*}
B_{n,1}(t) & = & \frac{1}{\sqrt{n}} \sum_{j=1}^n \textrm{e}^{t^\top X_j} \left(S_n^{-1/2}-\textrm{I}_d\right) X_j,\\
B_{n,2}(t) & = & \frac{1}{\sqrt{n}} \sum_{j=1}^n \textrm{e}^{t^\top X_j} \left(S_n^{-1/2}-\textrm{I}_d\right) \overline{X}_n, \ \
B_{n,3}(t)  =  \frac{1}{\sqrt{n}} \sum_{j=1}^n \textrm{e}^{t^\top X_j} \overline{X}_n.
\end{eqnarray*}
Since
\[
\|B_{n,2}(t)\|^2 \le \Big{(}\frac{1}{n}\sum_{j=1}^n \textrm{e}^{t^\top X_j} \Big{)}^2 \|\sqrt{n}\left(S_n^{-1/2}-\textrm{I}_d\right)\|_2^2 \cdot \|\overline{X}_n\|^2
\]
and $\|\sqrt{n}\left(S_n^{-1/2}-\textrm{I}_d\right)\|_2^2 \cdot \|\overline{X}_n\|^2 = O_\mathbb{P}(n^{-1})$, one may use (\ref{eabsch}) and Fubini's theorem to show
$B_{n,2} = o_\mathbb{P}(1)$. To proceed, rewrite $B_{n,3}(t)$ in the form
\[
B_{n,3}(t) = \frac{1}{n}\sum_{j=1}^n \textrm{e}^{t^\top X_j} \cdot \frac{1}{\sqrt{n}} \sum_{j=1}^n X_j.
\]
Since replacing the first factor with its expectation $\exp(\|t\|^2/2)$ means adding a term that is asymptotically negligible, we have
\begin{equation}\label{repbn3}
B_{n,3}(t) = \textrm{e}^{\|t\|^2/2} \, \frac{1}{\sqrt{n}} \sum_{j=1}^n X_j + o_\mathbb{P}(1).
\end{equation}
To tackle $B_{n,1}(t)$, we rewrite its transpose in the form
\[
B_{n,1}(t)^\top = \frac{1}{n}\sum_{j=1}^n \textrm{e}^{t^\top X_j} X_j^\top \, \sqrt{n}\left(S_n^{-1/2}-\textrm{I}_d\right)
\]
and use display (2.13) of \cite{HW1997}, according to which
\begin{equation}\label{darstsnx}
\sqrt{n} \left(S_n^{-1/2}-\textrm{I}_d\right) = - \frac{1}{2\sqrt{n}} \sum_{j=1}^n \left(X_jX_j^\top - \textrm{I}_d \right) + D_n,
\end{equation}
where $D_n = O_\mathbb{P}(n^{-1/2})$. Putting
\begin{equation}\label{defet}
E(t) =  \mathbb{E}\left[\textrm{e}^{t^\top X_1} X_1\right] = \frac{\textrm{d}}{\textrm{d}t} \, \textrm{e}^{\|t\|^2/2} = \textrm{e}^{\|t\|^2/2} \cdot t
\end{equation}
and $Y_n(t) = n^{-1} \sum_{j=1}^n \left(\exp(t^\top X_j) X_j - E(t) \right)$,
we have
\begin{equation}\label{repbn1}
B_{n,1}(t)^\top = \left(Y_n(t)  + E(t)\right)^\top \bigg{\{}- \frac{1}{2\sqrt{n}} \sum_{j=1}^n \left(X_jX_j^\top - \textrm{I}_d \right) + D_n\bigg{\}}.
\end{equation}
Abbreviating the matrix within curly brackets by $S_n$, we have $\|S_n Y_n(t)\|^2 \le \|S_n\|^2_2 \, \|Y_n(t)\|^2$. Since $\|S_n\|^2_2 = O_\mathbb{P}(1)$, it follows that
\[
\int_{\mathbb{R}^d} \|S_n Y_n(t)\|^2 \, w_\gamma(t) \, \textrm{d} t \ \le \ O_\mathbb{P}(1) \int_{\mathbb{R}^d} \|Y_n(t)\|^2 \, w_\gamma(t) \, \textrm{d} t.
\]
Now, observe that $Y_n(t)$ is a mean of centred random vectors, and invoking Fubini's theorem the expectation of the integral is seen to be of order $O(n^{-1})$. Thus,
$\|S_nY_n\|_\oplus = o_\mathbb{P}(1)$, and hence (since the matrix $D_n$ figuring in (\ref{repbn1}) is asymptotically negligible) we have
\[
B_{n,1}(t) = - \frac{1}{2\sqrt{n}} \sum_{j=1}^n \left(X_jX_j^\top - \textrm{I}_d \right) E(t) +o_\mathbb{P}(1).
\]
Upon combing this result with (\ref{repbn3}) and (\ref{defet}), the assertion follows.
\bewend

\begin{prop}\label{propant} For $A_n$ given in (\ref{ant}), we have
\[
A_n(t) = - \frac{1}{\sqrt{n}} \frac{\textrm{e}^{\|t\|^2/2}}{2} \sum_{j=1}^n \left(X_jX_j^\top t- t \right) + o_\mathbb{P}(1).
\]
\end{prop}

\noindent {\sc Proof.} Notice that $A_n(t) = A_{n,1}(t)-A_{n,2}(t)- A_{n,3}(t)$, where
\begin{eqnarray*}
A_{n,1}(t) & = & \frac{1}{\sqrt{n}} \sum_{j=1}^n \textrm{e}^{t^\top X_j} t^\top \left(S_n^{-1/2}- \textrm{I}_d \right) X_j(X_j-t),\\
A_{n,2}(t) & = & \frac{1}{\sqrt{n}} \sum_{j=1}^n \textrm{e}^{t^\top X_j} t^\top \left(S_n^{-1/2}- \textrm{I}_d \right) \overline{X}_n(X_j-t),\\
A_{n,3}(t) & = & \frac{1}{\sqrt{n}} \sum_{j=1}^n \textrm{e}^{t^\top X_j} t^\top  \overline{X}_n (X_j-t).
\end{eqnarray*}
To show $A_{n,2} = o_\mathbb{P}(1)$, use $n \big{\|}S_n^{-1/2}-\textrm{I}_d\big{\|}^2_2 \, \|\overline{X}_n\|^2 = O_\mathbb{P}(n^{-1})$ together with
\[
\|A_{n,2}(t)\|^2 \le \Big{\|}\frac{1}{n} \sum_{j=1}^n \textrm{e}^{t^\top X_j}(X_j-t)\Big{\|}^2 \cdot \|t\|^2 \cdot n \big{\|}S_n^{-1/2}-\textrm{I}_d\big{\|}^2_2 \cdot
\|\overline{X}_n\|^2
\]
and Fubini's theorem, since
$\mathbb{E} \big{\|}n^{-1} \sum_{j=1}^n \exp(t^\top X_j)(X_j-t)\big{\|}^2  = O(n^{-1})$,
due to the fact that the summands are centred random vectors. Likewise,
\[
\|A_{n,3}(t)\|^2 = \frac{1}{n^3} \sum_{i,j,k,\ell=1}^n \textrm{e}^{t^\top X_i} (X_i-t)^\top \cdot \textrm{e}^{t^\top X_j}(X_j-t) \cdot t^\top X_k \cdot t^\top X_\ell.
\]
Since each of the summands is a product of centred random  vectors or random variables, we have $\mathbb{E} \|A_{n,3}(t)\|^2 = O(1/n)$, and Fubini's theorem
yields $A_{n,3} = o_\mathbb{P}(1)$. To conclude the proof, observe that, by (\ref{darstsnx}),
\begin{equation}\label{an1darst}
A_{n,1}(t) = - \frac{1}{2\sqrt{n}} \sum_{j=1}^n \Delta_n(t) V_j(t) + \Delta_n(t) D_n \, t,
\end{equation}
where
$
\Delta_n(t) = n^{-1} \sum_{i=1}^n \exp(t^\top X_i) (X_i-t)X_i^\top$, $V_j(t) = X_jX_j^\top t - t$.
Notice that  $\Delta_n(t)$ is a mean of i.i.d. random matrices, and that
 \[
 \mathbb{E} \left[ \textrm{e}^{t^\top X_1} (X_1-t)X_1^\top  \right] = \textrm{e}^{\|t\|^2} \, \textrm{I}_d.
 \]
Straightforward calculations show that replacing $\Delta_n(t)$ with the right-hand side of the last equation means adding a term that is asymptotically negligible. Hence, the first term on the right-hand side of (\ref{an1darst}) is
\[
- \frac{1}{\sqrt{n}} \frac{\textrm{e}^{\|t\|^2/2}}{2} \sum_{j=1}^n \left(X_jX_j^\top t- t \right) + o_\mathbb{P}(1).
\]
The second term is $o_\mathbb{P}(1)$ since $D_n = O_\mathbb{P}(n^{-1/2})$.
\bewend

\begin{prop}{\rm{(}}Calculation of $K(s,t)${\rm{)}} \label{calckernel}\\
Recall $\widetilde{Z}_1(t)$ from (\ref{widetildez}).  Putting $m(s) = \exp(\|s\|^2/2)$ and $X_1 = X$, we have
\begin{eqnarray*}
\widetilde{Z}_1(s) \widetilde{Z}_1(t)^\top \! & \! = \! & \!
\textrm{e}^{(s+t)^\top X}(X-s)(X-t)^\top - m(s) \textrm{e}^{t^\top X} (X X^\top s-s)(X-t)^\top \\
\! & \!  \! &   - m(s) \textrm{e}^{t^\top X} X (X-t)^\top - m(t) \textrm{e}^{s^\top X} (X-s) (t^\top X X^\top - t^\top )\\
\! & \!  \! &   + m(s)m(t) \, (X X^\top s\! -\! s) (t^\top X X^\top \! - \! t^\top)  \\
\! & \!  \! &   + m(s)m(t) X  \, (t^\top X X^\top \! - \! t^\top) - m(t) \textrm{e}^{s^\top X} (X\! -\! s)X^\top \! \\
\! & \!  \! &   + m(s)m(t) (XX^\top \! s\! -\! s)X^\top \! + m(s)m(t) X X^\top.
\end{eqnarray*}
\end{prop}
Writing $\mathbf{0}$ for the zero matrix of  order $d$, we have
\begin{eqnarray*}
\mathbb{E}\left[ \textrm{e}^{(s+t)^\top X}(X-s)(X-t)^\top \right] & = & \textrm{e}^{\|s+t\|^2/2} \left(\textrm{I}_d + ts^\top \right),\\
\mathbb{E} \left[ \textrm{e}^{t^\top X} (X X^\top s-s)(X-t)^\top \right] & = & \textrm{e}^{\|t\|^2/2} \left(t s^\top + s^\top t \, \textrm{I}_d \right) ,\\
\mathbb{E} \left[ \textrm{e}^{t^\top X} X (X-t)^\top \right] & = & \textrm{e}^{\|t\|^2/2} \textrm{I}_d ,\\
\mathbb{E} \left[ \textrm{e}^{s^\top X} (X-s) (t^\top X X^\top - t^\top ) \right] & = & \textrm{e}^{\|s\|^2/2} \left(t s^\top + s^\top t \, \textrm{I}_d \right),\\
\mathbb{E} \left[ (X X^\top s\! -\! s) (t^\top X X^\top \! - \! t^\top) \right] & = & t s^\top + s^\top t \, \textrm{I}_d,\\
\mathbb{E} \left[ X  \, (t^\top X X^\top \! - \! t^\top) \right] & = & \mathbf{0},\\
\mathbb{E} \left[ \textrm{e}^{s^\top X} (X\! -\! s)X^\top \right] & = &  \textrm{e}^{\|s\|^2/2} \, \textrm{I}_d,\\
\mathbb{E} \left[ (XX^\top \! s\! -\! s)X^\top \right] & = & \mathbf{0},\\
\mathbb{E} \left[ X X^\top \right] & = & \textrm{I}_d.
\end{eqnarray*}
The assertion now follows from straightforward calculations.
\bewend

\begin{prop}\label{proplinalgebra} Let $(A_n)$ be a sequence of symmetric positive definite
($d \times d$)-matrices and $(b_n)$ an increasing sequence of positive real numbers satisfying
$\lim_{n\to \infty} b_n = \infty$ so that
\[
\lim_{n\to \infty} b_n \|A_n - {\rm{ I}}_d\|_2 = 0.
\]
We then have $\lim_{n\to \infty} b_n \|A_n^{-1/2} - {\rm{I}}_d\|_2 = 0$.
\end{prop}

\noindent {\sc Proof.} Let $\Lambda_n = \textrm{diag}(\lambda_1,\ldots,\lambda_n)$ be the diagonal matrix consisting of the
positive eigenvalues of $A_n$ so that $\|A_n\|_2 = \max_{i=1,\ldots,n} \lambda_i$ and
\begin{equation}\label{eigenungl}
\|(A_n^{1/2} + \textrm{I}_d)^{-1}\|_2 = \max_{i=1,\ldots,n} \left(\lambda_i^{1/2} +1\right)^{-1} < 1.
\end{equation}
Since the assumptions imply $\|A_n - {\rm{I}}_d\|_2 \to 0$,
choose $n_0$ so large that $\|A_n - {\rm{I}}_d\|_2 \le 1/2$ for each $n \ge n_0$.
Putting $T_n = \textrm{I}_d-A_n$, we have
\[
\|A_n^{-1}\|_2 = \|(\textrm{I}_d-T_n)^{-1}\|_2 \le \frac{1}{1-\|T_n\|_2} \le \frac{1}{1- 1/2} = 2, \quad n \ge n_0,
\]
and thus $\|A_n^{-1/2}\|_2 \le  \sqrt{2}, n\ge n_0$. Now, $A_n- \textrm{I}_d = (A_n^{1/2}+ \textrm{I}_d)A_n^{1/2} (\textrm{I}_d-A_n^{-1/2})$ implies
$\textrm{I}_d - A_n^{-1/2} = A_n^{-1/2}(A_n^{1/2}+\textrm{I}_d)^{-1}(A_n- \textrm{I}_d)$, whence
\[
\|A_n^{-1/2} -\textrm{I}_d\|_2 = \|\textrm{I}_d-A_n^{-1/2}\|_2 \le \|A_n^{-1/2}\|_2 \cdot \|(A_n^{1/2}+ \textrm{I}_d)^{-1}\|_2 \cdot \|A_n- \textrm{I}_d\|_2.
\]
In view of (\ref{eigenungl}) and $\|A_n^{-1/2}\|_2 \le  \sqrt{2}, n\ge n_0$, the assertion follows.
\bewend

% BibTeX users please use one of
%\bibliographystyle{spbasic}      % basic style, author-year citations
%\bibliography{}   % name your BibTeX data base

% Non-BibTeX users please use

\end{document}